\documentclass[12pt, A4paper]{amsart}
  
\usepackage{latexsym}
\usepackage[all]{xy}
\usepackage{amsthm}
\usepackage{amssymb}
\usepackage{color}
\usepackage{hyperref}
\usepackage{lscape}

\mathchardef\mhyphen="2D
\def\Mnd{\mathsf{Mnd}}
\def\2-Cat{2\mhyphen\mathsf{Cat}}
\def\cK{\mathcal K}
\def\ocK{\overline \cK}
\def\ox{\otimes}
\def\Wdl{\mathsf{Wdl}}
\def\two{\mathbf 2}

\newtheorem{proposition}{Proposition}[section]
\newtheorem*{proposition*}{Proposition}
\newtheorem*{lemma*}{Lemma}
\newtheorem{lemma}[proposition]{Lemma}
\newtheorem{corollary}[proposition]{Corollary}
\newtheorem{theorem}[proposition]{Theorem}

\theoremstyle{definition}
\newtheorem{definition}[proposition]{Definition}

\newtheorem{aussage}[proposition]{}

\numberwithin{equation}{section}

\newcounter{c}

\newcommand{\etyk}[1]{\vspace{-7.4mm}$$\begin{equation}\Label{#1}
\addtocounter{c}{1}}
\renewcommand{\]}{\ifnum \value{c}=1 $$\else \end{equation}\fi}
\begin{document}

\title{On the iteration of weak wreath products}

\author{Gabriella B\"ohm}
\address{Research Institute for Particle and Nuclear Physics, Budapest,
\newline\indent H-1525
Budapest 114, P.O.B.\ 49, Hungary}
\email{G.Bohm@rmki.kfki.hu}
\date{v1:Oct. 2011. v2:Jan 2012.}

\begin{abstract}
Based on a study of the 2-category of weak distributive laws, we describe a
method of iterating Street's weak wreath product construction in
\cite{Str:weak_dl}. That is, for any 2-category $\cK$ and for any non-negative
integer $n$, we introduce 2-categories $\Wdl^{(n)}(\cK)$, of $(n+1)$-tuples of
monads in $\cK$ pairwise related by weak distributive laws obeying the
Yang-Baxter equation. The first instance $\Wdl^{(0)}(\cK)$ coincides with
$\Mnd(\cK)$, the usual 2-category of monads in $\cK$, and for other values of
$n$, $\Wdl^{(n)}(\cK)$ contains $\Mnd^{n+1}(\cK)$ as a full 2-subcategory. For
the local idempotent closure $\ocK$ of $\cK$, extending the multiplication of
the 2-monad $\Mnd$, we equip these 2-categories with $n$ possible `weak wreath
product' 2-functors $\Wdl^{(n)}(\ocK)\to \Wdl^{(n-1)}(\ocK)$, such that all of
their possible $n$-fold composites $\Wdl^{(n)}(\ocK)\to \Wdl^{(0)}(\ocK)$ are
equal; i.e. such that the weak wreath product is `associative'.    
Whenever idempotent 2-cells in $\cK$ split, this leads to pseudofunctors
$\Wdl^{(n)}(\cK)\to \Wdl^{(n-1)}(\cK)$ obeying the associativity property
up-to isomorphism. 
We present a practically important occurrence of an iterated weak wreath
product: the algebra of observable quantities in an Ising type quantum spin
chain where the spins take their values in a dual pair of finite weak Hopf
algebras. 
We also construct a fully faithful embedding of $\Wdl^{(n)}(\ocK)$ into the
2-category of commutative $n+1$ dimensional cubes in $\Mnd(\ocK)$ (hence into
the 2-category of commutative $n+1$ dimensional cubes in $\cK$ whenever $\cK$
has Eilenberg-Moore objects and its idempotent 2-cells split).
Finally, we give a sufficient and necessary condition on a monad in $\ocK$ to
be isomorphic to an $n$-ary weak wreath product. 
\end{abstract}
\maketitle

\section*{Introduction}

At the heart of the iteration of wreath products in the work
\cite{Cheng:dlaws} of Eugenia Cheng, lies the 2-monad $\Mnd$ on the 2-category
$\2-Cat$ of 2-categories, 2-functors and 2-natural transformations, first
discussed in \cite{street}. For any 2-category $\cK$, the iteration of its
associative multiplication $\Mnd^{n}(\cK)\to \Mnd^{n-1}(\cK) \to \dots \to
\Mnd(\cK)$ takes an $(n+1)$-tuple of monads, pairwise related by distributive
laws obeying the Yang-Baxter equality, to a unique monad in $\cK$. The
resulting monad can be interpreted as an iterated wreath product. 

The aim of this paper is to find a similar iteration process for weak wreath
products introduced by Ross Street in \cite{Str:weak_dl} and by Stefaan
Caenepeel and Erwin De Groot in \cite{Stef&Erwin}.

These weak wreath products are defined in 2-categories in which idempotent
2-cells split, see \cite{Str:weak_dl}.
They are induced by weak distributive laws; i.e. certain 2-cells relating two
monads. They obey the usual compatibility conditions of distributive laws with
the multiplications of the monads, but the compatibility conditions with the
units are weakened \cite{Stef&Erwin}, \cite{Str:weak_dl}. 
Making weak distributive laws conceptually different from their non-weak
counterparts, they are not known to be monads in any 2-category. 
(However, a weak distributive law can be characterized as a pair of monads in
2-categories extending $\Mnd(\cK)$ and its variant $\Mnd^\ast(\cK)$,
respectively, see \cite{Boh_wtm}). 

In Section \ref{sec:iterate}, for any 2-category $\cK$, we construct a
2-category $\Wdl^{(n)}(\cK)$ for every non-negative integer $n$. Its objects
are $(n+1)$-tuples of monads in $\cK$ pairwise related by weak distributive
laws obeying the Yang-Baxter equation. The first one, $\Wdl^{(0)}(\cK)$ is
isomorphic to $\Mnd(\cK)$, the 2-category of monads in $\cK$ as defined in
\cite{street}. The next one, $\Wdl^{(1)}(\cK)$ is the 2-category of weak
distributive laws, obtained by dualizing the definition in
\cite{BoLaSt_wdlaw}. For every $n$, $\Wdl^{(n)}(\cK)$ contains
$\Mnd^{n+1}(\cK)$ as a full 2-subcategory. But, in contrast to the classical
(i.e. non-weak) case, $\Wdl^{(n)}(\cK)$ is not known to arise by the
$(n+1)$-fold application of some 2-monad. Although in this way we can not
interpret them as multiplications of some 2-monad, for each value of $n$ we
describe $n$ different 2-functors $\Wdl^{(n)}(\ocK) \to \Wdl^{(n-1)}(\ocK)$
(where $\ocK$ denotes the local idempotent closure of $\cK$). They extend the
$n$ possible multiplications $\Mnd^{n+1}(\ocK)\to \Mnd^n(\ocK)$. They give
rise to a unique composite $\Wdl^{(n)}(\ocK) \to \Wdl^{(0)}(\ocK)$ whose value
on an object of $\Wdl^{(n)}(\ocK)$ is regarded as the (associatively iterated)
weak wreath product of the $n+1$ occurring monads in $\ocK$. Whenever
idempotent 2-cells in $\cK$ split; that is, $\ocK$ and $\cK$ are biequivalent,
our construction yields pseudofunctors $\Wdl^{(n)}(\cK) \to \Wdl^{(n-1)}(\cK)$
giving rise to a composite $\Wdl^{(n)}(\cK) \to \Wdl^{(0)}(\cK)$ which is
unique up-to a pseudonatural equivalence in the choice of the biequivalence
$\ocK \to \cK$.   

Our motivation to study iterated weak wreath products comes from mathematical
physics. The Ising model is a quantum spin chain in which the spins take their
values in the sign group $\mathbb{Z}(2)$. In its various generalizations, the
spins may take their values in arbitrary finite groups \cite{SzlVec:G-spin},
in finite dimensional Hopf algebras \cite{NilSzl:Hopf_spin} or in finite
dimensional weak Hopf algebras \cite{NilSzlWie:Jones}, \cite{Bohm:PhD}. In all
of these, except the last quoted family of models, the algebra of the
observable quantities in any finite interval is given by an iterated wreath
product. In quantum spin chains based on weak Hopf algebras, however, the
algebras of observables are iterated weak wreath products. In Section
\ref{sec:spin} we present this example in some detail. 

The definition of $\Wdl^{(n)}(\ocK)$ is further motivated in Section
\ref{sec:embed} by a fully faithful embedding of it into the 2-category of
$n+1$ dimensional cubes in the 2-category of monads in $\ocK$. 
Whenever idempotent 2-cells in $\cK$ split, this yields a fully faithful
embedding of $\Wdl^{(n)}(\cK)$ into the 2-category of $n+1$ dimensional cubes
in the 2-category of monads in $\cK$; and also into the 2-category of $n+1$
dimensional cubes in $\cK$ whenever in addition $\cK$ admits Eilenberg-Moore
objects for monads.  

In our final Section \ref{sec:factor} we analyze the $n$-ary factorization
problem associated to weak distributive laws. That is, we give a complete
characterization of those monads in the local idempotent closure of a
2-category which arise as iterated weak wreath products of $n$ monads.   

Throughout, for a technical simplification, we work with 2-categories. There
is no difficulty, however, to extend our considerations to bicategories. 
\bigskip

\noindent
{\bf Acknowledgement.}
I would like to express my gratitude to Ross Street for his helpful comments
on this work. 
Partial support by the Hungarian Scientific Research Fund OTKA
K68195 is gratefully acknowledged.

\section{Preliminaries on weak distributive laws}
\label{sec:WDL}

In this section we revisit some recent `weak' generalizations of the formal
theory of monads that will be used in the sequel.

\begin{aussage}{\bf Local idempotent closure.}
To any 2-category $\cK$ we associate another 2-category $\ocK$ by freely
splitting idempotent 2-cells. In more detail, the 0-cells of $\ocK$ are the
same as those in $\cK$. The 1-cells in $\ocK$ are pairs consisting of a 1-cell
$v$ and a 2-cell $\overline v:v\to v$ in $\cK$ such that $\overline
v.\overline v=\overline v$; i.e. $\overline v$ is idempotent. The 2-cells
$(v,\overline v)\to (v',\overline v')$ in $\ocK$ are 2-cells $\omega:v\to v'$
in $\cK$ such that $\overline v'.\omega=\omega=\omega.\overline v$. Horizontal
and vertical compositions in $\ocK$ are induced by those in $\cK$. The
identity 2-cell is $\overline v:(v,\overline v)\to (v,\overline v)$.

Throughout, we shall use the notation seen above: If it is not otherwise
stated, in a 1-cell in $\ocK$, for the idempotent 2-cell part we use the
overlined version of the same symbol that denotes the 1-cell part.

For any 2-category $\cK$, there is an evident inclusion 2-functor $\cK\to
\ocK$, acting on the 0-cells as the identity map, taking a 1-cell $v$ to
$(v,v)$ -- i.e. the 1-cell with identity 2-cell part -- and acting on the
2-cells again as the identity map. 

We say that idempotent 2-cells in a 2-category $\cK$ split if, for any
idempotent 2-cell $\theta: v\to v$ there exist a 1-cell $w$ and 2-cells
$\iota:w\to v$ and $\pi:v\to w$ such that $\pi.\iota=w$ and
$\iota.\pi=\theta$. If the splitting exists then it is unique up-to
an isomorphism of $w$. Clearly, in $\ocK$ idempotent 2-cells split for any
2-category $\cK$.

Whenever in $\cK$ idempotent 2-cells split, the inclusion $\cK\to \ocK$
becomes a biequivalence. (Since it acts on the 0-cells as the identity map,
this simply means that it induces an equivalence of the hom categories.) Hence
there is a {\em pseudofunctor} $\ocK\to \cK$ which is its inverse (in the
sense of inverse biequivalences). 
On the 0-cells also this pseudofunctor acts as the identity map. On a
1-cell $(v,\overline v)$ its action is constructed via a {\em chosen}
splitting of the idempotent 2-cell $\overline v$. If $(\iota:w\to v,\pi:v\to
w)$ is this chosen splitting, then the image of $(v,\overline v)$ is $w$. A
2-cell $\omega:(v,\overline v)\to (v',\overline {v'})$ is taken to 
$
\xymatrix @C=20pt{
w\ar[r]|(.45){\, \iota\,}&
v\ar[r]|(.45){\,\omega\,}&
v'\ar[r]|(.45){\raisebox{4pt}{${}_{\,\pi'\,}$}}&
w'}
$.
Let us stress that the biequivalence $\ocK\to \cK$ is not a 2-functor in
general and it is unique only up-to a pseudonatural equivalence arising from
the choice of the splitting of each idempotent 2-cell. 
\end{aussage}

\begin{aussage}{\bf Demimonads.} \label{as:demimonad}
In simplest terms, a demimonad in a 2-category is a monad $(A,(t,\overline
t))$ in the local idempotent closure, cf. \cite{BoLaSt_weak}. Explicitly, it
is given by a 1-cell $t:A\to A$ and 2-cells $\mu:t^2 \to t$ and $\eta:1_A\to
t$ such that the following diagrams commute. 
$$
\xymatrix{
t^3\ar[r]^-{\mu t}\ar[d]_-{t\mu}&
t^2\ar[d]^-\mu\\
t^2\ar[r]_-\mu&
t
}
\quad 
\xymatrix{
t\ar[r]^-{\eta t}\ar[d]_-{t\eta}&
t^2\ar[d]^-\mu\\
t^2\ar[r]_-\mu&
t
}\quad
\xymatrix{
1_A \ar[r]^-{\eta \eta}\ar[dr]_-{\eta}&
t^2\ar[d]^-\mu\\
&t
}\quad
\xymatrix{
t^2\ar[r]^-{\eta t^2}\ar[rd]_-{\mu}&
t^3\ar[d]^-{\mu^2}\\
&t
} 
$$
By the unitality condition, the idempotent 2-cell $\overline t$ must be equal
to $\mu.t\eta=\mu.\eta t$ (hence it is a redundant information that will be
often omitted in the sequel). This structure occurred in \cite{Boh_wtm} under
the name `pre-monad'.  

A demimonad $(A,(t,\overline t))$ is the image of a monad under the inclusion
$\cK\to \ocK$ if and only if $\overline t$ is the identity 2-cell $t$.
\end{aussage}

\begin{aussage}{\bf Weak distributive laws.} \label{as:wdl}
Extending the notion of distributive law due to Jon Beck (see \cite{Beck}),
weak distributive laws in a 2-category  were introduced by Ross Street in
\cite{Str:weak_dl} as follows.  
They consist of two monads $(A,t)$ and $(A,s)$ on the same object, and a
2-cell $\lambda:ts\to st$ such that the following diagrams 
commute. 
\begin{equation}\label{eq:wdl}
\xymatrix@C=12pt{
t^2s\ar[r]^-{t\lambda}\ar[d]^-{\mu s}&
tst\ar[r]^-{\lambda t}&
st^2\ar[d]_-{s\mu}\\
ts\ar[rr]_--\lambda&&
st
}\ 
\xymatrix@C=12pt{
s\ar[rrr]^-{\eta s}\ar[d]^-{s \eta}&&&
ts\ar[d]_-\lambda\\
st\ar[r]_-{st \eta}&
sts\ar[r]_-{s\lambda}&
s^2t\ar[r]_-{\mu t}&
st
}\ 
\xymatrix@C=12pt{
ts^2\ar[r]^-{\lambda s}\ar[d]^-{t\mu}&
sts\ar[r]^-{s\lambda}&
s^2t\ar[d]_-{\mu t}\\
ts\ar[rr]_--\lambda&&
st
}\ 
\xymatrix@C=12pt{
t\ar[rrr]^-{t\eta}\ar[d]^-{\eta t}&&&
ts\ar[d]_-\lambda\\
st\ar[r]_-{\eta st}&
tst\ar[r]_-{\lambda t}&
st^2\ar[r]_-{s\mu}&
st
}
\end{equation}
The same set of axioms occurred also in \cite{Stef&Erwin}.
By \cite[Proposition 2.2]{Str:weak_dl}, the second and fourth diagrams can
be replaced by a single diagram
\begin{equation}\label{eq:wunit}
\xymatrix{
st\ar[r]^-{\eta st}\ar[d]_-{st\eta}&
tst\ar[r]^-{\lambda t}&
st^2\ar[d]^-{s\mu}\\
sts\ar[r]_-{s\lambda}&
s^2t\ar[r]_-{\mu t}&
st\ .
}
\end{equation}
The equal paths around \eqref{eq:wunit} give rise to an idempotent 2-cell
$\overline{\lambda}:st\to st$ (which occurs also in the bottom rows of the 
second and fourth diagrams in \eqref{eq:wdl}). It is an identity if and only
if $\lambda$ is a distributive law in the strict sense. 

Note that a weak distributive law in $\cK$ is the same as a weak distributive
law in the horizontal opposite of $\cK$. 

A weak distributive law in $\ocK$ is then given by demimonads $(A,t)$ and
$(A,s)$ and a 2-cell $\lambda: ts\to st$ rendering commutative the diagrams in
\eqref{eq:wdl} and obeying in addition the normalization conditions
$$
\xymatrix @C=30pt {
ts\ar[r]^-{t \mu.t\eta s}\ar[d]_-\lambda\ar[rd]^-\lambda&
ts\ar[d]^-\lambda&&
ts\ar[r]^-{\mu s.t\eta s}\ar[d]_-\lambda\ar[rd]^-\lambda&
ts\ar[d]^-\lambda\\
st\ar[r]_-{\mu t.s\eta t}&
st&&
st\ar[r]_-{s\mu.s\eta t}&
st\ .}
$$
In the sequel we shall need some identities on weak distributive laws (in
$\ocK$).  
The axioms imply commutativity of the following diagrams, see 
\cite{Str:weak_dl}. 
\begin{equation}\label{eq:wdlid}
\xymatrix @R=10pt {
ts\ar[r]^-{\lambda}\ar[rdd]_-{\lambda}&
st\ar[dd]^-{\overline \lambda}
&
s^2t^2\ar[r]^-{s\overline{\lambda}t}\ar[dd]_-{\mu\mu}&
s^2t^2\ar[dd]^-{\mu\mu}
&
sts\ar[r]^-{s\lambda}\ar[d]_-{\overline \lambda s}&
s^2t\ar[dd]^-{\mu t}
&
tst\ar[r]^-{\lambda t}\ar[d]_-{t\overline\lambda}&
st^2\ar[dd]^-{s\mu}\\
&&&&
sts\ar[d]_-{s\lambda}&
&
tst\ar[d]_-{\lambda t}\\
&st
&
st\ar[r]_-{\overline \lambda}&
st
&
s^2t\ar[r]_-{\mu t}&
st
&
st^2\ar[r]_-{s\mu}&
st}
\end{equation}
Moreover, by the associativity of $\mu$, the left-bottom path in the last
diagram in (\ref{eq:wdl}) commutes with the multiplication by $t$ on 
the right. Hence so does the top-right path meaning the commutativity of the
first diagram in  
\begin{equation}\label{eq:star}
\xymatrix{
t^2\ar[d]_-\mu\ar[r]^-{t\eta t}&
tst\ar[r]^-{\lambda t}&
st^2\ar[d]^-{s\mu}&&
s^2\ar[d]_-\mu\ar[r]^-{s\eta s}&
sts\ar[r]^-{s\lambda }&
s^2t\ar[d]^-{\mu t}\\
t\ar[r]_-{t\eta}&
ts\ar[r]_-\lambda&
st&&
s\ar[r]_-{\eta s}&
ts\ar[r]_-\lambda&
st.}
\end{equation}
Commutativity of the second diagram follows symmetrically.
\end{aussage}

\begin{aussage}{\bf The 2-category of weak distributive laws.} 
\label{as:2-Cat_wdl}
Dualizing in the appropriate sense the definition of the 2-category of {\em
mixed} weak distributive laws in \cite{BoLaSt_wdlaw}, the following 2-category
$\Wdl(\cK)$ of weak distributive laws in $\cK$ is obtained (see
\cite[Paragraph 1.9]{BohmGomTor}). The 0-cells are the weak distributive laws
$\lambda:ts\to st$. The 1-cells $\lambda\to \lambda'$ are triples consisting
of a 1-cell $v:A\to A'$ and 2-cells $\xi:t'v\to vt$ and $\zeta:s'v\to vs$ in
$\cK$, such that $(v,\xi):(A,t)\to (A',t')$ and $(v,\zeta):(A,s)\to (A',s')$
are  1-cells in $\Mnd(\cK)$ (also called {\em monad morphisms} in
\cite{street}) and the following diagram commutes.  
\begin{equation}\label{eq:Wdl_1_cell}
\xymatrix @C=6pt{
t's'v\ar[rrr]^-{t'\zeta}\ar[d]_-{\lambda' v}&&&
t'vs\ar[rrr]^-{\xi s}&&&
vts\ar[d]^-{v\lambda}\\
s't'v\ar[rr]_-{s'\xi}&&
s'vt\ar[rr]_-{\zeta t}&&
vst\ar[rr]_-{v\overline{\lambda}}&&
vst
}
\end{equation}
The 2-cells $(v,\xi,\zeta)\to(v',\xi',\zeta')$ are 2-cells $\omega: v\to v'$
in $\cK$ which are 2-cells in $\Mnd(\cK)$ (i.e. {\em monad transformations} by
the terminology of \cite{street}); both $(v,\xi)\to (v',\xi')$ and
$(v,\zeta)\to (v',\zeta')$. 
Horizontal and vertical compositions are induced by those in $\cK$.
This definition can be interpreted in terms of (weak) liftings as in
\cite{BoLaSt_wdlaw}. 

There is a fully faithful embedding $\Mnd^2(\cK)\to \Wdl(\cK)$ as
follows. It takes a 0-cell $((A,t),(s,\lambda))$ to the distributive law
$\lambda:ts\to st$, regarded as a weak distributive law. It takes a 1-cell
$((v,\xi),\zeta)$ to $(v,\xi,\zeta)$ and it acts on the 2-cells as the
identity map. 
\end{aussage}

\begin{aussage}{\bf Weak wreath product.} \label{as:w_wreath}
The weak wreath product induced by a weak distributive law 
in a 2-category in which idempotent 2-cells split, 
was discussed by Ross Street in \cite[Theorem 2.4]{Str:weak_dl}. In the
particular case of the monoidal category (i.e. one object bicategory) of
modules over a commutative ring, it appeared in \cite[Theorem
  3.2]{Stef&Erwin}. 

For an arbitrary 2-category $\cK$, there is a {\em weak wreath product}
2-functor $\Wdl(\ocK)\to \Mnd(\ocK)$, which sends a weak distributive law
$\lambda:ts\to st$ to the monad $(st,\overline{\lambda})$ in $\overline{\cK}$,
with multiplication and unit 
$$
\xymatrix{
(st)^2\ar[r]^-{s\lambda t}&
s^2t^2\ar[r]^-{\mu\mu}&
st}
\quad \textrm{and}\quad
\xymatrix{
1\ar[r]^-{\eta\eta}&
ts\ar[r]^-{\lambda}&
st\ .}
$$
It sends a 1-cell $((v,\overline v),\xi,\zeta):\lambda\to \lambda'$ to the
monad morphism with the same 1-cell part $(v,\overline v)$ and the 2-cell part 
$$
\xymatrix{
s't'v\ar[r]^-{s'\xi}&
s'vt\ar[r]^-{\zeta t}&
vst\ar[r]^-{v\overline{\lambda}}&
vst\ .
}
$$
On the 2-cells it acts as the identity map.

Whenever idempotent 2-cells in $\cK$ split, the biequivalence $\ocK\simeq
\cK$ induces a pseudo\-functor $\Wdl(\cK)\stackrel\simeq \to \Wdl(\ocK)\to
\Mnd(\ocK) \stackrel\simeq \to\Mnd(\cK)$. (It can be chosen, in fact, to be a
2-functor by choosing the biequivalence $\ocK\to \cK$ adopting the convention
that we split identity 2-cells trivially; i.e. via identity 2-cells.) Its
object map yields Street's weak wreath product in $\cK$.
\end{aussage}

\begin{aussage}{\bf Binary factorization.} \label{as:bin_wfac}
Let $\cK$ be any 2-category. As proved in 
\cite{BohmGomTor}, a demimonad $(A,r)$ is isomorphic to a weak wreath product
induced by some weak distributive law $ts\to st$ in $\ocK$ if and only if the
following hold. 
\begin{itemize}
\item[{(a)}] There are 1-cells in $\Mnd(\ocK)$ with trivial 1-cell parts
$$
\xymatrix{
(A,(t,\overline t))&&
(A,(r,\overline r))\ar[ll]_-{((A,A),\alpha)} \ar[rr]^-{((A,A),\beta)}&&
(A,(s,\overline s))\ ;}
$$
\item[{(b)}] The 2-cell
$$
\pi:=\big(
\xymatrix{
(st,\overline s \overline t)\ar[r]^-{\beta\alpha}&
(rr,\overline r \,\overline r)\ar[r]^-\mu&
(r,\overline r)}\big)
$$
in $\ocK$ possesses section $\iota$ (meaning $\pi.\iota=\overline r\equiv
\mu.r\eta$) which is an $s$-$t$ bimodule morphism with respect to the $t$- and
$s$-actions induced on $r$ by $\alpha$ and $\beta$, respectively. 
\end{itemize}
Indeed, for the weak wreath product induced by a weak distributive law
$\lambda:ts\to st$, we have 1-cells 
$$
\xymatrix{
(A,(t,\overline t))&&
(A,(st,\overline\lambda))
\ar[ll]_-{((A,A),\lambda.t\eta)}\ar[rr]^-{((A,A),\lambda.\eta s)}&&
(A,(s,\overline s))}
$$
in $\Mnd(\ocK)$. Moreover, the 2-cell $\pi$ in part (b) comes out as
$$
\big(\xymatrix{
(st,\overline s \overline t)\ar[rr]^-{\lambda\lambda.\eta st\eta}&&
(stst,\overline \lambda\,\overline \lambda)\ar[rr]^-{\mu\mu.s\lambda t}&&
(st,\overline \lambda)}\big)=
\big(\xymatrix{
(st,\overline s\overline t)\ar[r]^-{\overline \lambda}&
(st,\overline \lambda)}\big)\ ,
$$
which is split by the bimodule morphism $\overline \lambda:(st,\overline
\lambda)\to (st,\overline s \overline t)$. 

Conversely, if properties (a) and (b) hold, then 
$$
\xymatrix{
ts\ar[r]^-{\alpha\beta}&
rr\ar[r]^-\mu&
r\ar[r]^-\iota&
st}
$$
is a weak distributive law with corresponding idempotent equal to
$\iota.\pi:st\to st$. The isomorphism between the induced weak wreath product
and $(A,r)$ is provided by $\xymatrix{
(st,\iota.\pi)\ar[r]<2pt>^-\pi&(r,\overline r)\ar[l]<2pt>^-\iota} $ in $\ocK$.
For the details of the proof we refer to \cite{BohmGomTor}.
\end{aussage}

\section{2-categories of weak distributive laws and the iterated weak wreath
  product} 
\label{sec:iterate}

Throughout this section, $\cK$ is an arbitrary 2-category and $\ocK$ stands
for its local idempotent closure. For any non-negative integer $n$, we define
a 2-category $\Wdl^{(n)}(\cK)$. Its objects are $(n+1)$-tuples of monads
pairwise related by weak distributive laws obeying the Yang-Baxter condition. 
For each value of $n$, we construct $n$
different 2-functors $\Wdl^{(n)}(\ocK)\to \Wdl^{(n-1)}(\ocK)$ corresponding to
taking the weak wreath product of two consecutive monads of the $n+1$
occurring ones. 
We show that these 2-functors give rise to a unique composite
$\Wdl^{(n)}(\ocK)\to \Wdl^{(0)}(\ocK)=\Mnd(\ocK)$. We regard its object map as
the $n$-ary weak wreath product of the involved monads.

\begin{aussage} {\bf The 2-category $\Wdl^{(n)}(\cK)$}.\label{as:Wdln}
For any non-negative integer $n$, a 0-cell of $\Wdl^{(n)}(\cK)$ 
is given by $n+1$ monads $s_0,s_1,\dots,s_n$ together with weak distributive
laws $\lambda_{i,j}:s_j s_i \to s_i s_j$ for all $0\leq i< j\leq n$, obeying
for all $0\leq i<j<k\leq n$ the Yang-Baxter relation 
$$
\xymatrix{
s_ks_js_i\ar[r]^-{\lambda_{j,k}s_i}\ar[d]_-{s_k \lambda_{i,j}}&
s_js_ks_i\ar[r]^-{s_j\lambda_{i,k}}&
s_js_is_k\ar[d]^-{\lambda_{i,j}s_k}\\
s_ks_is_j\ar[r]_-{\lambda_{i,k}s_j}&
s_is_ks_j\ar[r]_-{s_i \lambda_{j,k}}&
s_is_js_k.}
$$ 
The 1-cells consist of a 1-cell $v$ and 2-cells $\xi_i:s'_iv\to vs_i$ in
$\cK$ for all $0\leq i\leq n$, such that $(v,\xi_i,\xi_j)$ is a 1-cell
$\lambda_{i,j}\to \lambda'_{i,j}$ in $\Wdl(\cK)$ (see Paragraph
\ref{as:2-Cat_wdl}), for all $0\leq i<j\leq n$.      
The 2-cells are those 2-cells $\omega:v\to v'$ in $\cK$ which are 2-cells
$(v,\xi_i,\xi_j) \to (v',\xi'_i,\xi'_j)$ in $\Wdl(\cK)$ (in the sense of
Paragraph \ref{as:2-Cat_wdl}), for all $0\leq i<j\leq n$.   
Since $\Wdl(\cK)$ is closed under the horizontal and vertical compositions in
$\cK$, so is $\Wdl^{(n)}(\cK)$. Hence it is a 2-category with the horizontal
and vertical compositions induced by those in $\cK$.

Recall from \cite{Cheng:dlaws} that a 0-cell in $\Mnd^n(\cK)$ is given by $n$
monads, pairwise related by distributive laws obeying the Yang-Baxter
condition. The 1-cells consist of $n$ monad morphisms for the $n$ involved
monads with a common underlying 1-cell, obeying \eqref{eq:Wdl_1_cell} (in the
simplified form when the occurring idempotents are identities). The 2-cells
are those 2-cells in $\cK$ which are monad transformations for all of the $n$
monad morphisms. With this description in mind, extending that in Paragraph
\ref{as:2-Cat_wdl}, there is an evident fully faithful embedding
$\Mnd^{n+1}(\cK)\to \Wdl^{(n)}(\cK)$. 
\end{aussage}

Taking any $m+1$-element subset of $\{0,1,\dots, n\}$ induces
an evident 2-functor $\Wdl^{(n)}(\cK)$ $\to \Wdl^{(m)}(\cK)$. 

\begin{lemma}
Take any object $\{\lambda_{i,j}:s_js_i \to s_is_j\}_{0\leq i<j\leq 2}$ of
$\Wdl^{(2)}(\ocK)$. 
In addition to $\overline s_i:s_i\to s_i$ for $0\leq i\leq 2$, and $\overline 
\lambda_{i,j}:s_is_j\to s_is_j$ for $0\leq i<j\leq 2$, let us introduce the
the following idempotent 2-cells in $\cK$. 
$$
\xymatrix @R=10pt{
{\overrightarrow {\lambda}}_{0,p,q}:=\big(
s_0s_ps_q\ar[r]^-{s_0s_ps_q\eta_0}&
s_0s_ps_qs_0\ar[r]^-{s_0s_p \lambda_{0,q}}&
s_0s_ps_0s_q\ar[r]^-{s_0\lambda_{0,p}s_q}&
s_0^2s_ps_q\ar[r]^-{\mu_0s_ps_q}&
s_0s_ps_q\big),\\
{\overleftarrow {\lambda}}_{k,l,2}
:=\big(
s_ks_ls_2\ar[r]^-{\eta_2s_ks_ls_2}&
s_2s_ks_ls_2\ar[r]^-{\lambda_{k,2}s_ls_2}&
s_ks_2s_ls_2\ar[r]^-{s_k\lambda_{l,2}s_2}&
s_ks_ls_2^2\ar[r]^-{s_ks_l\mu_2}&
s_ks_ls_2\big),}
$$
for $p=1,q=2$ and $p=2,q=1$; and for $k=0,l=1$ and $k=1,l=0$. 
They obey the following equalities.
\begin{eqnarray}
&&\overrightarrow \lambda_{0,p,q}.\overline \lambda_{0,p}s_q=
\overrightarrow \lambda_{0,p,q}=
\overline \lambda_{0,p}s_q.\overrightarrow \lambda_{0,p,q}
\label{eq:rarrownorm}\\
&&\overrightarrow \lambda_{0,1,2}.s_0\lambda_{1,2}=
s_0\lambda_{1,2}.\overrightarrow \lambda_{0,2,1}
\label{eq:rarrowintw}\\
&&\overrightarrow \lambda_{0,p,q}.\lambda_{0,p}s_q=
\lambda_{0,p}s_q.s_p\overline \lambda_{0,q}
\label{eq:rarrow_bar}\\
&&\overleftarrow \lambda_{k,l,2}.s_k\overline \lambda_{l,2}=
\overleftarrow \lambda_{k,l,2}=
s_k\overline \lambda_{l,2}.\overleftarrow \lambda_{k,l,2}
\label{eq:larrownorm}\\
&&\overleftarrow \lambda_{0,1,2}.\lambda_{0,1}s_2=
\lambda_{0,1}s_2.\overleftarrow \lambda_{1,0,2}
\label{eq:larrowintw}\\
&&\overleftarrow \lambda_{k,l,2}.s_k\lambda_{l,2}=
s_k\lambda_{l,2}.\overline \lambda_{k,2}s_l
\label{eq:larrow_bar}\\
&&\overleftarrow \lambda_{0,1,2}.\overrightarrow \lambda_{0,1,2}=
\overleftarrow \lambda_{0,1,2}.\overline \lambda_{0,1}s_2=
\overrightarrow \lambda_{0,1,2}.s_0\overline \lambda_{1,2}=
\overrightarrow \lambda_{0,1,2}.\overleftarrow \lambda_{0,1,2}
\label{eq:3_idemp}
\ .
\end{eqnarray}
\end{lemma}

In what follows, we shall denote by $\overline \lambda_{012}$ the equal
2-cells in \eqref{eq:3_idemp} (we shall see later the irrelevance of inserting 
any comma between the labels). 

\begin{proof}
We only present a proof of \eqref{eq:3_idemp}, verification of the other
equalities is left to the reader.  

The first and the last expressions in \eqref{eq:3_idemp} are equal by
commutativity of the following diagram. 
$$
\xymatrix @C=30pt @R=15pt{
s_0s_1s_2\ar[r]^-{s_0s_1s_2\eta_0}\ar[d]_-{\eta_2 s_0s_1s_2}&
s_0s_1s_2s_0\ar[r]^-{s_0s_1\lambda_{0,2}}&
s_0s_1s_0s_2\ar[r]^-{s_0\lambda_{0,1}s_2}&
s_0^2s_1s_2 \ar[r]^-{\mu_0 s_1s_2}&
s_0s_1s_2\ar[d]^-{\eta_2 s_0s_1s_2}\\
s_2s_0s_1s_2\ar[r]^-{s_2s_0s_1s_2\eta_0}\ar[d]_-{\lambda_{0,2}s_1s_2}&
s_2s_0s_1s_2s_0\ar[r]^-{s_2s_0s_1\lambda_{0,2}}&
s_2s_0s_1s_0s_2\ar[r]^-{s_2s_0\lambda_{0,1}s_2}&
s_2s_0^2s_1s_2 \ar[r]^-{s_2\mu_0 s_1s_2}\ar[d]^-{\lambda_{0,2}s_0s_1s_2}
\ar@{}[rdd]|-{\eqref{eq:wdl}}&
s_2s_0s_1s_2\ar[dd]^-{\lambda_{0,2}s_1s_2}\\
s_0s_2s_1s_2\ar[r]^-{s_0s_2 s_1s_2\eta_0}\ar[d]_-{s_0\lambda_{1,2}s_2}&
s_0s_2s_1s_2s_0\ar[r]^-{s_0s_2 s_1\lambda_{0,2}}&
s_0s_2s_1s_0s_2\ar[r]^-{s_0s_2
  \lambda_{0,1}s_2}\ar[d]^-{s_0\lambda_{1,2}s_0s_2}
\ar@{}[rdd]|-{\mathrm{(YB)}}&
s_0s_2 s_0s_1s_2 \ar[d]^-{s_0\lambda_{0,2}s_1s_2}\\
s_0s_1s_2^2\ar[r]^-{s_0s_1s_2^2\eta_0}\ar[dd]_-{s_0s_1\mu_2}&
s_0s_1s_2^2s_0\ar[r]^-{s_0s_1s_2\lambda_{0,2}}\ar[dd]_-{s_0s_1\mu_2s_0}
\ar@{}[rdd]|-{\eqref{eq:wdl}}&
s_0s_1s_2s_0s_2\ar[d]^-{s_0s_1\lambda_{0,2}s_2}&
s_0^2s_2s_1s_2\ar[r]^-{\mu_0s_2s_1s_2}\ar[d]^-{s_0^2\lambda_{1,2}s_2}&
s_0s_2s_1s_2\ar[d]^-{s_0\lambda_{1,2}s_2}\\
&&s_0s_1s_0s_2^2\ar[r]^-{s_0\lambda_{0,1}s_2^2}\ar[d]^-{s_0s_1s_0\mu_2}&
s_0^2s_1s_2^2\ar[r]^-{\mu_0s_1s_2^2}\ar[d]^-{s_0^2s_1\mu_2}&
s_0s_1s_2^2\ar[d]^-{s_0s_1\mu_2}\\
s_0s_1s_2\ar[r]_-{s_0s_1s_2\eta_0}&
s_0s_1s_2s_0\ar[r]_-{s_0s_1\lambda_{0,2}}&
s_0s_1s_0s_2\ar[r]_-{s_0\lambda_{0,1}s_2}&
s_0^2s_1s_2 \ar[r]_-{\mu_0 s_1s_2}&
s_0s_1s_2}
$$
The first and the third expressions in \eqref{eq:3_idemp} are equal by
commutativity of the following diagram.
$$
\xymatrix @C=30pt @R=15pt{
s_0s_1s_2\ar[r]^-{s_0s_1s_2\eta_0}\ar[d]_-{s_0\eta_2 s_1s_2}&
s_0s_1s_2s_0\ar[r]^-{s_0s_1\lambda_{0,2}}&
s_0s_1s_0s_2\ar[r]^-{s_0\lambda_{0,1}s_2}&
s_0^2s_1s_2 \ar[r]^-{\mu_0 s_1s_2}\ar[d]^-{s_0\eta_2 s_0s_1s_2}
\ar@{}[rdd]|-{\eqref{eq:star}}&
s_0s_1s_2\ar[d]^-{\eta_2 s_0s_1s_2}\\
s_0s_2s_1s_2\ar[r]^-{s_0s_2s_1s_2\eta_0}\ar[d]_-{s_0\lambda_{1,2}s_2}&
s_0s_2s_1s_2s_0\ar[r]^-{s_0s_2s_1\lambda_{0,2}}&
s_0s_2s_1s_0s_2\ar[r]^-{s_0s_2\lambda_{0,1}s_2}\ar[d]^-{s_0\lambda_{1,2}s_0s_2}
\ar@{}[rdd]|-{\mathrm{(YB)}}&
s_0s_2s_0s_1s_2 \ar[d]^-{s_0\lambda_{0,2}s_1s_2}&
s_2s_0s_1s_2 \ar[d]^-{\lambda_{0,2}s_1s_2}\\
s_0s_1s_2^2\ar[r]^-{s_0s_1s_2^2\eta_0}\ar[dd]_-{s_0s_1\mu_2}&
s_0s_1s_2^2s_0\ar[r]^-{s_0s_1s_2\lambda_{0,2}}\ar[dd]_-{s_0s_1\mu_2s_0}
\ar@{}[rdd]|-{\eqref{eq:wdl}}&
s_0s_1s_2s_0s_2\ar[d]^-{s_0s_1\lambda_{0,2}s_2}&
s_0^2s_2s_1s_2\ar[r]^-{\mu_0s_2s_1s_2}\ar[d]^-{s_0^2\lambda_{1,2}s_2}&
s_0s_2s_1s_2\ar[d]^-{s_0\lambda_{1,2}s_2}\\
&&s_0s_1s_0s_2^2\ar[r]^-{s_0\lambda_{0,1}s_2^2}\ar[d]^-{s_0s_1s_0\mu_2}&
s_0^2s_1s_2^2\ar[r]^-{\mu_0s_1s_2^2}&
s_0s_1s_2^2\ar[d]^-{s_0s_1\mu_2}\\
s_0s_1s_2\ar[r]_-{s_0s_1s_2\eta_0}&
s_0s_1s_2s_0\ar[r]_-{s_0s_1\lambda_{0,2}}&
s_0s_1s_0s_2\ar[r]_-{s_0\lambda_{0,1}s_2}&
s_0^2s_1s_2 \ar[r]_-{\mu_0 s_1s_2}&
s_0s_1s_2}
$$
Equality of the second and last expressions in \eqref{eq:3_idemp} follows
symmetrically. 
\end{proof}

\begin{lemma}\label{lem:C_0}
For any object $\{\lambda_{i,j}:s_js_i \to s_is_j\}_{0\leq i<j\leq 2}$ of
$\Wdl^{(2)}(\ocK)$, consider the monads $(s_0s_1,\overline \lambda_{0,1})$ and
$(s_1s_2,\overline\lambda_{1,2})$ in $\ocK$, induced by the weak distributive
laws $\lambda_{0,1}$ and $\lambda_{1,2}$, respectively. There are weak
distributive laws  
$$
\xymatrix @R=10pt{
\lambda_{01,2}:=\big(
s_2 (s_0 s_1) \ar[r]^-{\lambda_{0,2} s_1}&
s_0 s_2 s_1 \ar[r]^-{s_0 \lambda_{1,2}}&
s_0 s_1 s_2 \ar[r]^-{\overrightarrow{\lambda}_{0,1,2}}
&
(s_0 s_1) s_2 \big)\\
{\color{white}\lambda_{01,2}\ }
=\big(s_2 (s_0 s_1) \ar[r]^-{s_2\overline\lambda_{01}}&
s_2s_0s_1\ar[r]^-{\lambda_{0,2} s_1}&
s_0 s_2 s_1 \ar[r]^-{s_0 \lambda_{1,2}}&
(s_0s_1)s_2\big)&
\textrm{and}\\
\lambda_{0,12}:=\big((s_1 s_2)s_0 \ar[r]^-{s_1 \lambda_{0,2}}&
s_1 s_0 s_2 \ar[r]^-{\lambda_{0,1} s_2}&
s_0 s_1 s_2\ar[r]^-{\overleftarrow{\lambda}_{0,1,2}}
&
s_0 (s_1 s_2)\big)\\ 
{\color{white} \lambda_{0,12}\ }
=\big((s_1 s_2)s_0 \ar[r]^-{\overline\lambda_{12}s_0}&
s_1s_2s_0\ar[r]^-{s_1 \lambda_{0,2}}&
s_1 s_0 s_2 \ar[r]^-{\lambda_{0,1} s_2}&
s_0 (s_1 s_2)\big)}
$$  
in $\ocK$. Moreover, their induced monads
$(s_0s_1s_2,\overline\lambda_{01,2})$ and
$(s_0s_1s_2,\overline\lambda_{0,12})$ are equal. \end{lemma}  
 
\begin{proof}
Both given forms of $\lambda_{01,2}$ are equal by \eqref{eq:rarrowintw} and
\eqref{eq:rarrow_bar}. It is a 2-cell in $\ocK$ by \eqref{eq:rarrownorm}. 
Compatibility with the multiplication of $s_2$ holds since both
$\lambda_{0,2}$ and $\lambda_{1,2}$ are compatible with it. 
With the normalization conditions $\overline \lambda_{01}. \mu_{01}=\mu_{01}=
\mu_{01}.\overline\lambda_{01}s_0s_1= \mu_{01}.s_0s_1\overline\lambda_{01}$ at
hand, compatibility of $\lambda_{01,2}$ with the multiplication of $s_0s_1$
follows by the compatibilities of $\lambda_{0,2}$ with $\mu_0$ and of
$\lambda_{1,2}$ with $\mu_1$ and the Yang-Baxter condition. 
The weak unitality condition \eqref{eq:wunit} follows by the equality of the
second and third expressions in \eqref{eq:3_idemp} (so that
$\overline\lambda_{01,2}= \overline\lambda_{012}$).
This proves that $\lambda_{01,2}$ is a weak distributive law and
$\lambda_{0,12}$ can be handled symmetrically (in particular,
$\overline\lambda_{0,12}= \overline\lambda_{012}$).   

Equality of the units in the induced monads follows immediately
by the Yang-Baxter condition. Concerning the multiplications, composing the
equal paths around 
$$
\xymatrix @R=15pt @C=22pt{
s_1s_2s_1s_0s_1\ar[rr]^-{s_1s_2\lambda_{0,1}s_1}
\ar[ddd]_-{s_1\lambda_{1,2}s_0s_1}\ar@{}[rrddd]|-{\mathrm{(YB)}}&&
s_1s_2s_0s_1^2\ar[rrr]^-{s_1s_2s_0\mu_1}\ar[dd]^-{s_1\lambda_{0,2}s_1^2}&&&
s_1s_2s_0s_1 \ar[d]^-{s_1\lambda_{0,2}s_1}\\
&&&&&
s_1 s_0 s_2 s_1 \ar[d]^-{\lambda_{0,1}s_2s_1}\\
&&
s_1s_0s_2s_1^2\ar[r]^-{\lambda_{0,1}s_2s_1^2}\ar[d]^(.4){s_1s_0\lambda_{1,2}s_1}&
s_0s_1s_2s_1^2\ar[rr]^-{s_0s_1s_2\mu_1}\ar[d]^(.4){s_0s_1\lambda_{1,2}s_1}&&
s_0s_1 s_2 s_1 \ar[d]^(.4){s_0s_1 \lambda_{1,2}}\\
s_1^2s_2s_0s_1\ar[r]^-{s_1^2\lambda_{0,2}s_1}\ar[d]_-{\mu_1s_2s_0s_1}&
s_1^2 s_0 s_2 s_1 \ar[r]^-{s_1 \lambda_{0,1}s_2s_1}
\ar[d]^-{\mu_1 s_0 s_2 s_1}&
s_1s_0s_1s_2s_1 \ar[r]^-{\lambda_{0,1}s_1s_2s_1}&
s_0s_1^2s_2s_1\ar[r]^-{s_0s_1^2\lambda_{1,2}}\ar[d]^-{s_0\mu_1s_2s_1}&
s_0s_1^3s_2\ar[r]^-{s_0s_1\mu_1s_2}\ar[d]^-{s_0\mu_1s_1s_2}&
s_0s_1^2s_2\ar[d]^-{s_0\mu_1s_2}\\
s_1s_2s_0s_1\ar[r]_-{s_1\lambda_{0,2}s_1}&
s_1s_0s_2s_1 \ar[rr]_-{\lambda_{0,1}s_2s_1}&&
s_0s_1s_2s_1\ar[r]_-{s_0s_1\lambda_{1,2}}&
s_0s_1^2s_2\ar[r]_-{s_0\mu_1s_2}&
s_0s_1s_2}
$$
by $s_1s_2\eta_1s_0s_1$ on the right, we obtain
\begin{equation}\label{eq:equal_mp}
s_0\mu_1s_2.\lambda_{0,1}\lambda_{1,2}.s_1\lambda_{0,2}s_1.
\overline{\lambda}_{12}s_0s_1=
s_0\mu_1s_2.\lambda_{0,1}\lambda_{1,2}.s_1\lambda_{0,2}s_1.
s_1s_2 \overline{\lambda}_{01}.
\end{equation}
Inserting these equal 2-cells $s_1s_2s_0s_1 \to s_0s_1s_2$ into
$\mu_0s_1\mu_2.s_0(-) s_2$, we conclude the equality of the
multiplications induced on $s_0s_1s_2$ by $\lambda_{01,2}$ and
$\lambda_{0,12}$, respectively. 
\end{proof}

\begin{aussage} {\bf On the Yang-Baxter condition.}
Actually, also a sort of converse of Lemma \ref{lem:C_0} holds. Consider weak
distributive laws $\{\lambda_{i,j}:s_js_i \to s_is_j\}_{0\leq i<j\leq 2}$ in
$\ocK$. Assume that 
$$
\xymatrix @R=10pt{
\lambda_{01,2}:=\big(
s_2 (s_0 s_1) \ar[r]^-{s_2\overline\lambda_{01}}&
s_2s_0s_1\ar[r]^-{\lambda_{0,2} s_1}&
s_0 s_2 s_1 \ar[r]^-{s_0 \lambda_{1,2}}&
(s_0s_1)s_2\big)&
\textrm{and}\\
\lambda_{0,12}:=\big(
(s_1 s_2)s_0 \ar[r]^-{\overline\lambda_{12}s_0}&
s_1s_2s_0\ar[r]^-{s_1 \lambda_{0,2}}&
s_1 s_0 s_2 \ar[r]^-{\lambda_{0,1} s_2}&
s_0 (s_1 s_2)\big)}
$$  
are weak distributive laws inducing equal monads
$(s_0s_1s_2,\overline\lambda_{01,2})$ and
$(s_0s_1s_2,\overline\lambda_{0,12})$. Then the Yang-Baxter condition holds. 

Indeed, equality of the multiplications $\mu_{0,12}$ and $\mu_{01,2}$ is
equivalent to \eqref{eq:equal_mp}. With this identity at hand, from the
compatibility of $\lambda_{0,12}$ with $\mu_{12}$ we obtain 
$$
\lambda_{0,12}.\mu_{12}s_0=
s_0\mu_1\mu_2.
\lambda_{0,1}s_1s_2^2.
s_1s_0\lambda_{1,2}s_2.
s_1\lambda_{0,2}s_1s_2.
s_1s_2\lambda_{0,12}\ .
$$
Precomposing this equality with $\eta_1s_2s_1\eta_2s_0$, we conclude that 
\begin{equation}\label{eq:YBL}
\lambda_{0,1}s_2.s_1\lambda_{0,2}.\lambda_{1,2}s_0=
s_0s_1\mu_2.
s_0\lambda_{1,2}s_2.
\lambda_{0,2}s_1s_2.
s_2\lambda_{0,1}s_2.
s_2s_1\lambda_{0,2}.
s_2\lambda_{1,2}s_0.
s_2\eta_2s_1s_0\ .
\end{equation}
Symmetrically,
\begin{equation}\label{eq:YBR}
s_0\lambda_{1,2}.\lambda_{0,2}s_1.s_2\lambda_{0,1}=
\mu_0s_1s_2.
s_0\lambda_{0,1}s_2.
s_0s_1\lambda_{0,2}.
s_0\lambda_{1,2}s_0.
\lambda_{0,2}s_1s_0.
s_2\lambda_{0,1}s_0.
s_2s_1\eta_0s_0\ .
\end{equation}
It follows by the associativity of $\mu_1$ that 
$s_0\overline \lambda_{1,2}.\overline\lambda_{0,1}s_2=
\overline\lambda_{0,1}s_2.s_0\overline \lambda_{1,2}$. Precomposing with
$\eta_0s_1\eta_2$, we obtain from this
$
s_0\lambda_{1,2}.
s_0\eta_2s_1.
\lambda_{0,1}.
s_1\eta_0=
\lambda_{0,1}s_2.
s_1\eta_0s_2.
\lambda_{1,2}.
\eta_2s_1$.
Inserting these latter equal expressions into
$$
s_0s_1\mu_2.
s_0\lambda_{1,2}s_2.
\lambda_{0,2}s_1s_2.
s_2\mu_0s_1s_2.
s_2s_0\lambda_{0,1}s_2.
s_2s_0s_1\lambda_{0,2}.
s_2(-)s_0
$$
and using \eqref{eq:YBL} and \eqref{eq:YBR} to simplify both sides of the
resulting equality, we obtain the Yang-Baxter condition.
\end{aussage}

\begin{lemma}\label{lem:C_1}
For any 1-cell $\{\xi_i:s'_iv\to vs_i\}_{0\leq i\leq 2}$ in
$\Wdl^{(2)}(\ocK)$, the following yield 1-cells in $\Wdl(\ocK)$ between 
the 0-cells described in Lemma \ref{lem:C_0}.
\begin{eqnarray*}
&&\{\xi_{01}:=\big(
\xymatrix{s'_0s'_1v\ar[r]^-{s'_0\xi_1}&
s'_0vs_1\ar[r]^-{\xi_0s_1}&
vs_0s_1\ar[r]^-{v\overline\lambda_{01}}&
vs_0s_1}\big),\qquad
\xymatrix{
s'_2v\ar[r]^-{\xi_2}&
vs_2}\}\quad\textrm{and}\\
&&\{\xymatrix{
s'_0v\ar[r]^-{\xi_0}&
vs_0},\qquad
\xi_{12}:=\big(
\xymatrix{s'_1s'_2v\ar[r]^-{s'_1\xi_2}&
s'_1vs_2\ar[r]^-{\xi_1s_2}&
vs_1s_2\ar[r]^-{v\overline\lambda_{12}}&
vs_1s_2}\big)\}
\end{eqnarray*}
\end{lemma}

\begin{proof}
By Paragraph \ref{as:w_wreath}, $\xi_{01}$ and $\xi_{12}$ are 1-cells in
$\Mnd(\ocK)$. Moreover, $\xi_{01}$ and $\xi_2$ obey \eqref{eq:Wdl_1_cell} by
commutativity of the following diagram. 
$$
\xymatrix @R=15pt @C=30pt {
s'_2s'_0s'_1v\ar[dd]_-{s'_2\xi_{01}}\ar[r]^-{s'_2\overline{\lambda'}_{01}v}&
s'_2s'_0s'_1v\ar[rr]^-{\lambda'_{0,2}s'_1v}\ar[d]^-{s'_2s'_0\xi_1}&&
s'_0s'_2s'_1v\ar[rr]^-{s'_0\lambda'_{1,2}v}\ar[d]^-{s'_0s'_2\xi_1}
\ar@{}[rrdd]|-{\eqref{eq:Wdl_1_cell}}&&
s'_0s'_1s'_2v\ar[d]^-{s'_0s'_1\xi_2}\\
&
s'_2s'_0vs_1\ar[rr]^-{\lambda'_{0,2}vs_1}\ar[d]^-{s'_2\xi_0s_1}
\ar@{}[rrdd]|-{\eqref{eq:Wdl_1_cell}}&&
s'_0s'_2vs_1\ar[d]^-{s'_0\xi_2s_1}&&
s'_0s'_1vs_2\ar[d]^-{s'_0\xi_1s_2}\\
s'_2vs_0s_1\ar[dd]_-{\xi_2s_0s_1}&
s'_2vs_0s_1\ar[l]_-{s'_2v\overline\lambda_{01}}\ar[d]^-{\xi_2s_0s_1}&&
s'_0vs_2s_1\ar[r]^-{s'_0v\lambda_{1,2}}\ar[d]^-{\xi_0s_2s_1}&
s'_0vs_1s_2\ar[d]^-{\xi_0s_1s_2}&
s'_0vs_1s_2\ar[l]_-{s'_0v\overline\lambda_{12}}\ar[d]^-{\xi_0s_1s_2}\\
&
vs_2s_0s_1\ar[r]^-{v\lambda_{0,2}s_1}\ar[d]^-{vs_2\overline\lambda_{01}}
\ar@{}[rd]|-{\quad \eqref{eq:rarrow_bar}}&
vs_0s_2s_1 \ar[d]^-{v\overrightarrow \lambda_{0,2,1}}
\ar@{}[rd]|-{\quad \eqref{eq:rarrownorm}}&
vs_0s_2s_1\ar[r]^-{vs_0\lambda_{1,2}}\ar[d]^-{v\overrightarrow\lambda_{0,2,1}}
\ar[l]_-{v\overline\lambda_{02}s_1}\ar@{}[rd]|-{\quad \eqref{eq:rarrowintw}}&
vs_0s_1s_2\ar[d]^-{v\overrightarrow\lambda_{0,1,2}}
\ar@{}[rd]|-{\quad \eqref{eq:3_idemp}}&
vs_0s_1s_2\ar[l]_-{vs_0\overline\lambda_{12}}\ar[d]^-{v\overline\lambda_{012}}\\
vs_2s_0s_1\ar[r]_-{vs_2\overline\lambda_{01}}&
vs_2s_0s_1\ar[r]_-{v\lambda_{0,2}s_1}&
vs_0s_2s_1\ar@{=}[r]&
vs_0s_2s_1\ar[r]_-{vs_0\lambda_{1,2}}&
vs_0s_1s_2\ar@{=}[r]&
vs_0s_1s_2}
$$
The top-left region commutes since $\xi_{01}$ is a 2-cell in $\Mnd(\ocK)$ of
domain $(s'_0s'_1v,\overline{\lambda'}_{0,1}\overline v)$.
Also $\xi_0$ and $\xi_{12}$ obey \eqref{eq:Wdl_1_cell}, hence
constitute a 1-cell in $\Wdl(\ocK)$, by commutativity of the similar diagram
below. The bottom-left region commutes by the normalization of $\xi_{12}$. 
$$
\xymatrix @R=15pt @C=30pt {
s'_1s'_2s'_0v\ar[r]^-{\overline{\lambda'}_{12}s'_0v}\ar[d]_-{s'_1s'_2\xi_0}&
s'_1s'_2s'_0v\ar[rr]^-{s'_1\lambda'_{0,2}v}\ar[d]^-{s'_1s'_2\xi_0}
\ar@{}[rrdd]|-{\eqref{eq:Wdl_1_cell}}&&
s'_1s'_0s'_2v\ar[rr]^-{\lambda'_{0,1}s'_2v}\ar[d]^-{s'_1s'_0\xi_2}&&
s'_0s'_1s'_2v\ar[d]^-{s'_0s'_1\xi_2}\\
s'_1s'_2vs_0\ar[ddd]_-{\xi_{12}s_0}\ar[r]^-{\overline{\lambda'}_{12}vs_0}&
s'_1s'_2vs_0\ar[d]^-{s'_1\xi_2s_0}&&
s'_1s'_0vs_2\ar[rr]^-{\lambda'_{0,1}vs_2}\ar[d]^-{s'_1\xi_0s_2}
\ar@{}[rrdd]|-{\eqref{eq:Wdl_1_cell}}&&
s'_0s'_1vs_2\ar[d]^-{s'_0\xi_1s_2}\\
&
s'_1vs_2s_0\ar[r]^-{s'_1v\lambda_{0,2}}\ar[d]^-{\xi_1s_2s_0}&
s'_1vs_0s_2\ar[d]^-{\xi_1s_0s_2}&
s'_1vs_0s_2\ar[l]_-{s'_1v\overline\lambda_{02}}\ar[d]^-{\xi_1s_0s_2}&&
s'_0vs_1s_2\ar[d]^-{\xi_0s_1s_2}\\
&
vs_1s_2s_0\ar[r]^-{vs_1\lambda_{0,2}}\ar[d]^-{v\overline\lambda_{12}s_0}
\ar@{}[rd]|-{\quad \eqref{eq:larrow_bar}}&
vs_1s_0s_2\ar[d]^-{v\overleftarrow\lambda_{1,0,2}}
\ar@{}[rd]|-{\quad \eqref{eq:larrownorm}}&
vs_1s_0s_2\ar[r]^-{v\lambda_{0,1}s_2}\ar[d]^-{v\overleftarrow\lambda_{1,0,2}}
\ar[l]_-{vs_1\overline\lambda_{02}}\ar@{}[rd]|-{\quad\eqref{eq:larrowintw}}&
vs_0s_1s_2\ar[d]^-{v\overleftarrow\lambda_{0,1,2}}
\ar@{}[rd]|-{\quad \eqref{eq:3_idemp}}&
vs_0s_1s_2 \ar[l]_-{v\overline\lambda_{01}s_2}\ar[d]^-{v\overline\lambda_{012}}\\
vs_1s_2s_0\ar[r]_-{v\overline\lambda_{12}s_0}&
vs_1s_2s_0\ar[r]_-{vs_1\lambda_{0,2}}&
vs_1s_0s_2\ar@{=}[r]&
vs_1s_0s_2\ar[r]_-{v\lambda_{0,1}s_2}&
vs_0s_1s_2\ar@{=}[r]&
vs_0s_1s_2}
$$
\end{proof}

\begin{theorem}\label{thm:n_w-wreath}
For any 2-category $\cK$, and any positive integer $n$,
there are $n$ different 2-functors $C_k:\Wdl^{(n)}(\ocK) \to
\Wdl^{(n-1)}(\ocK)$, for $1\leq k\leq n$, as follows. They take a 0-cell
$\{\lambda_{i,j}:s_js_i\to s_is_j\}_{0\leq i<j\leq n}$ to  
\begin{equation}\label{eq:Ck0}
\left\{
\begin{array}{ll}
\xymatrix{
s_js_i\ar[r]^-{\lambda_{i,j}}&s_is_j}&
\mathrm{if}\ i,j\notin\{k-1,k\}\\
\xymatrix{s_j (s_{k-1} s_k) \ar[r]^-{s_j\overline\lambda_{k-1,k}}&
s_js_{k-1}s_k\ar[r]^-{\lambda_{k-1,j} s_k}&
s_{k-1} s_j s_k \ar[r]^-{s_{k-1} \lambda_{k,j}}&
(s_{k-1} s_k) s_j}&\mathrm{if}\ k<j\\
\xymatrix{
(s_{k-1} s_k)s_i \ar[r]^-{\overline\lambda_{k-1,k}s_i}&
s_{k-1}s_ks_i\ar[r]^-{s_{k-1} \lambda_{i,k}}&
s_{k-1} s_i s_k \ar[r]^-{\lambda_{i,k-1} s_k}&
s_i (s_{k-1} s_k)}&\mathrm{if}\ i<k-1
\end{array}
\right\}
\end{equation}
where $(s_{k-1} s_k,\overline\lambda_{k-1,k})$ is the monad in $\ocK$ induced
by the weak distributive law $\lambda_{k-1,k}$. They send a 1-cell $\{ \xi_i:
s'_iv\to vs_i\}_{0\leq i\leq n}$ to  
\begin{equation}\label{eq:Ck1}
\left\{
\begin{array}{ll}
\xymatrix{s'_iv\ar[r]^-{\xi_i}&vs_i}&
\mathrm{if}\ 0\leq i<k-1\\
\xymatrix{
(s'_{k-1}s'_k)v\ar[r]^-{s'_{k-1}\xi_k}&
s'_{k-1}vs_k\ar[r]^-{\xi_{k-1}s_k}&
vs_{k-1}s_k\ar[r]^-{v\overline \lambda_{k-1,k}}&
v(s_{k-1}s_k)}&\\
\xymatrix{s'_{i+1}v\ar[r]^-{\xi_{i+1}}&vs_{i+1}}&
\mathrm{if}\ k-1<i<n
\end{array}
\right\}
\end{equation}
On the 2-cells they act as the identity map.
\end{theorem}

\begin{proof}
By Lemma \ref{lem:C_0}, each line in \eqref{eq:Ck0} is a weak
distributive law in $\ocK$. We only need to check the Yang-Baxter
conditions. For $0\leq i<j<k-1<n$ the Yang-Baxter condition follows by
commutativity of 
$$
\xymatrix @C=45pt @R=18pt {
s_{k-1}s_ks_js_i\ar[r]^-{\overline\lambda_{k-1,k}s_js_i}
\ar[rd]^(.6){\overline\lambda_{k-1,k}s_js_i}
\ar[d]_-{s_{k-1}s_k\lambda_{i,j}}&
s_{k-1}s_ks_js_i \ar[r]^-{s_{k-1}\lambda_{j,k}s_i}
&
s_{k-1}s_js_ks_i\ar[r]^-{\lambda_{j,k-1}s_ks_i}&
s_js_{k-1}s_ks_i\ar[d]^-{s_j\overline \lambda_{k-1,k}s_i}\\
s_{k-1}s_ks_is_j\ar[rd]^(.6){\overline\lambda_{k-1,k}s_is_j}
\ar[d]_-{\overline\lambda_{k-1,k}s_is_j}
&
s_{k-1}s_ks_js_i\ar[r]^-{s_{k-1}\lambda_{j,k}s_i}\ar[d]^-{s_{k-1}s_k\lambda_{i,j}}
\ar@{}[rdd]|-{\mathrm{(YB)}}&
s_{k-1}s_js_ks_i
\ar[r]^-{\lambda_{j,k-1}s_ks_i}\ar[d]^-{s_{k-1}s_j\lambda_{i,k}}&
s_js_{k-1}s_ks_i\ar[d]^-{s_js_{k-1}\lambda_{i,k}}\\
s_{k-1}s_ks_is_j\ar[d]_-{s_{k-1}\lambda_{i,k}s_j}&
s_{k-1}s_ks_is_j\ar[d]^-{s_{k-1}\lambda_{i,k}s_j}&
s_{k-1}s_js_is_k\ar[r]^-{\lambda_{j,k-1}s_is_k}
\ar[d]^-{s_{k-1}\lambda_{i,j}s_k}\ar@{}[rdd]|-{\mathrm{(YB)}}&
s_js_{k-1}s_is_k\ar[d]^-{s_j\lambda_{i,k-1}s_k}\\
s_{k-1}s_is_ks_j\ar[d]_-{\lambda_{i,k-1} s_ks_j}&
s_{k-1}s_is_ks_j\ar[r]^-{s_{k-1}s_i\lambda_{j,k}}\ar[d]^-{\lambda_{i,k-1}s_ks_j}&
s_{k-1}s_is_js_k\ar[d]^-{\lambda_{i,k-1}s_js_k}&
s_js_is_{k-1}s_k\ar[d]^-{\lambda_{i,j}s_{k-1}s_k}\\
s_is_{k-1}s_ks_j\ar[r]_-{s_i\overline \lambda_{k-1,k}s_j}&
s_is_{k-1}s_ks_j\ar[r]_-{s_is_{k-1}\lambda_{j,k}}&
s_is_{k-1}s_js_k\ar[r]_-{s_i\lambda_{j,k-1}s_k}&
s_is_js_{k-1}s_k\ .}
$$
The top-right region and the bottom-left region commute by Lemma \ref{lem:C_0}.
The $0< k<i<j\leq n$ case is treated symmetrically. The Yang-Baxter condition
in the last case, when $0\leq i<k-1$ and $k<j\leq n$, follows by commutativity
of the similar diagram
$$
\xymatrix @C=28pt @R=18pt {
s_js_{k-1}s_ks_i\ar[r]^-{s_j\overline\lambda_{k-1,k}s_i}
\ar[rd]^(.6){s_j\overline\lambda_{k-1,k}s_i}
\ar[d]_-{s_j\overline\lambda_{k-1,k}s_i}&
s_js_{k-1}s_ks_i\ar[r]^-{\lambda_{k-1,j}s_ks_i}&
s_{k-1}s_js_ks_i\ar[r]^-{s_{k-1}\lambda_{k,j}s_i}&
s_{k-1}s_ks_js_i\ar[r]^-{s_{k-1}s_k\lambda_{i,j}}
\ar[d]^-{\overline \lambda_{k-1,k}s_js_i}&
s_{k-1}s_ks_is_j\ar[d]^-{\overline \lambda_{k-1,k}s_is_j}\\
s_js_{k-1}s_ks_i\ar[d]_-{s_js_{k-1}\lambda_{i,k}}&
s_js_{k-1}s_ks_i\ar[r]^-{\lambda_{k-1,j}s_ks_i}\ar[d]^-{s_js_{k-1}\lambda_{i,k}}&
s_{k-1}s_js_ks_i\ar[r]^-{s_{k-1}\lambda_{k,j}s_i}
\ar[d]^-{s_{k-1}s_j\lambda_{i,k}}\ar@{}[rrd]|-{\mathrm{(YB)}}& 
s_{k-1}s_ks_js_i  \ar[r]^-{s_{k-1}s_k\lambda_{i,j}}&
s_{k-1}s_ks_i s_j\ar[d]^-{s_{k-1}\lambda_{i,k}s_j}\\
s_js_{k-1}s_is_k\ar[d]_-{s_j\lambda_{i,k-1}s_k}&
s_js_{k-1}s_is_k\ar[r]^-{\lambda_{k-1,j}s_is_k}\ar[d]^-{s_j\lambda_{i,k-1}s_k}
\ar@{}[rrdd]|-{\mathrm{(YB)}}&
s_{k-1}s_js_is_k\ar[r]^-{s_{k-1}\lambda_{i,j}s_k}&
s_{k-1}s_is_js_k\ar[r]^-{s_{k-1}s_i\lambda_{k,j}}
\ar[dd]^-{\lambda_{i,k-1}s_js_k}&
s_{k-1}s_is_ks_j\ar[dd]^-{\lambda_{i,k-1}s_ks_j}\\
s_js_is_{k-1}s_k \ar[r]^-{s_js_i\overline{\lambda}_{k-1,k}}
\ar[d]_-{\lambda_{i,j}s_{k-1}s_k}&
s_js_is_{k-1}s_k\ar[d]^-{\lambda_{i,j}s_{k-1}s_k}&
&&\\
s_is_js_{k-1}s_k\ar[r]_-{s_is_j\overline \lambda_{k-1,k}}&
s_is_js_{k-1}s_k\ar[rr]_-{s_i\lambda_{k-1,j}s_k}&
&
s_is_{k-1}s_js_k\ar[r]_-{s_is_{k-1}\lambda_{k,j}}&
s_is_{k-1}s_ks_j\ .}
$$
Both regions at the top-left commute by Lemma \ref{lem:C_0}.
This proves that \eqref{eq:Ck0} describes a 0-cell in $\Wdl^{(n-1)}(\ocK)$.
By Lemma \ref{lem:C_1}, \eqref{eq:Ck1} is a 1-cell in $\Wdl^{(n-1)}(\ocK)$.
Evidently, 2-cells in $\Wdl^{(n)}(\ocK)$ are 2-cells in $\Wdl^{(n-1)}(\ocK)$
as well. Hence the stated maps define 2-functors $C_k$ which clearly preserve
the horizontal and vertical compositions.
\end{proof}

Via the fully faithful embedding $\Mnd^{n+1}(\ocK)\to \Wdl^{(n)}(\ocK)$ in
Paragraph \ref{as:Wdln}, the 2-functors in Theorem \ref{thm:n_w-wreath} extend
the multiplication $C$ of the 2-monad $\Mnd$. That is, the following diagram
commutes, for all $1\leq k\leq n$.
$$
\xymatrix @C=40pt{
\Mnd^{n+1}(\ocK)\ar[rr]^-{\Mnd^{k-1}C\Mnd^{n-k}(\ocK)}\ar[d]&&
\Mnd^{n}(\ocK)\ar[d]\\
\Wdl^{(n)}(\ocK)\ar[rr]_-{C_k}&&
\Wdl^{(n-1)}}
$$
By associativity of the 2-monad $\Mnd$, the $n$-fold iteration of the
2-functor in the top row; i.e.
$$
\xymatrix @C=40pt{ 
\Mnd^{n+1}(\ocK)
\ar[rr]^-{\Mnd^{k_n-1}C\Mnd^{n-k_n}(\ocK)}&&
\quad \dots \quad 
\ar[rr]^-{\Mnd^{k_2-1}C\Mnd^{2-k_2}(\ocK)}&&
\Mnd^2(\ocK)
\ar[r]^-C&
\Mnd(\ocK)}
$$
does not depend on the values of $k_i\in\{1,\dots,i\}$, for $1\leq i\leq
n$. That is, its object map describes an `associative' wreath product of
monads. Although the 2-functor in the bottom row is not known to correspond to
the multiplication in any 2-monad, in the rest of this section we show that it
describes an associative weak wreath product in an appropriate sense. 

\begin{lemma}\label{lem:n_idemp}
For any integer $n>1$, and for any 0-cell $\{\lambda_{i,j}:s_js_i \to
s_is_j\}_{0\leq i< j \leq n}$ of $\Wdl^{(n)}(\ocK)$, consider the idempotent
2-cell 
$$
\scalebox{.95}{
\xymatrix @C=22pt @R=6pt{
{\overleftarrow{\lambda}}_{0,1,\dots,n}:=\big(s_0s_1\dots s_n 
\ar[rr]^-{\eta_ns_0s_1\dots s_n}&&
s_ns_0s_1\dots s_n\ar[r]&
s_0s_1\dots s_{n-1}s_n^2\ar[rr]^-{s_0s_1\dots s_{n-1}\mu_n}&&
s_0s_1\dots s_n\big)}}
$$
in $\cK$ (where the unlabelled arrow denotes the unique composite of
$\lambda_{i,j}$s of the given domain and codomain) and
\begin{equation}\label{eq:n_idemp}
\overline\lambda_{01\dots n}:=
\overleftarrow \lambda_{0,1,\dots,n}.
\overleftarrow \lambda_{0,1,\dots,n-1}s_n.\ \cdots\ .
\overleftarrow \lambda_{0,1,2}s_3s_4\dots s_n.
\overline \lambda_{01}s_2s_3\dots s_n\ .
\end{equation}
This construction in \eqref{eq:n_idemp} associates the same idempotent 2-cell
to any 0-cell $\{\lambda_{i,j}:s_js_i \to s_is_j\}_{0\leq i< j \leq n}$ of
$\Wdl^{(n)}(\ocK)$ and to its image  under any of the 2-functors
$C_k$ in Theorem \ref{thm:n_w-wreath}. That is, for all $1\leq k\leq n$,
$\overline\lambda_{01\dots n}= \overline\lambda_{01\dots k-2(k-1,k)k+1\dots n}$.
\end{lemma}

\begin{proof}
Using commutativity of
\begin{equation}\label{eq:lambda_n_idemp}
\scalebox{.9}{
\xymatrix {
s_{k+1}s_0\dots s_k
\ar[r]_-{\raisebox{-8pt}{${}_{s_{k+1}\eta_ks_0\dots s_k}$}}
\ar@/^1.5pc/[rrrr]^-{s_{k+1}\overleftarrow\lambda_{0,\dots,k}}
\ar@/_2pc/[rddd]_-{s_{k+1}\overleftarrow\lambda_{0,\dots,k}}&
s_{k+1}s_ks_0\dots s_k\ar[rr]\ar[dd]\ar[rd]&&
s_{k+1}s_0\dots s_{k-1}s_k^2
\ar[r]_-{\raisebox{-8pt}{${}_{s_{k+1}s_0\dots s_{k-1}\mu_k}$}}\ar[d]
\ar@{}[rd]|-{\eqref{eq:wdl}}&
s_{k+1}s_0\dots s_k\ar[d]\\
&&
s_ks_0\dots s_{k+1}\ar[r]\ar[rd]
\ar@{}[ru]|-{\mathrm{(YB)}}&
s_0\dots s_{k-1}s_k^2s_{k+1}
\ar[r]_-{\raisebox{-5pt}{${}_{s_0\dots s_{k-1}\mu_ks_{k+1}}$}}
\ar[d]^-{\overleftarrow \lambda_{0,\dots,k}s_ks_{k+1}}&
s_0\dots s_{k+1}\ar[dd]_-{\overleftarrow \lambda_{0,\dots,k}s_{k+1}}\\
&
s_{k+1}s_0\dots s_{k-1}s_k^2\ar[rr]\ar[d]^-{s_{k+1}s_0\dots s_{k-1}\mu_k}
\ar@{}[ru]|-{\mathrm{(YB)}}\ar@{}[rrd]|-{\eqref{eq:wdl}}&&
s_0\dots s_{k-1}s_k^2s_{k+1}\ar[d]^-{s_0\dots s_{k-1}\mu_ks_{k+1}}\\
&
s_{k+1}s_0\dots s_k\ar[rr]&&
s_0\dots s_{k+1}\ar@{=}[r]&
s_0\dots s_{k+1}}}
\end{equation}
for any $1\leq k < n$, it follows easily that 
$\overline\lambda_{0\dots n}$ is idempotent. 
By \eqref{eq:Ck0} on one hand, and by \eqref{eq:lambda_n_idemp}
on the other,
\begin{eqnarray*}
&&\overleftarrow\lambda_{0,\dots,k-2,(k-1,k),k+1,\dots m}=
\overleftarrow\lambda_{0,\dots,m}.
s_0\dots s_{k-2}\overline\lambda_{k-1,k}s_{k+1}\dots s_m\\
&&\hspace{3.2cm} =s_0\dots s_{k-2}\overline\lambda_{k-1,k}s_{k+1}\dots s_m.
\overleftarrow\lambda_{0,\dots,m}.
s_0\dots s_{k-2}\overline\lambda_{k-1,k}s_{k+1}\dots s_m
\end{eqnarray*}
for all $k<m\leq n$. Moreover, by commutativity of 
$$
\xymatrix @C=40pt @R=15pt {
s_0\dots s_k\ar[r]^-{\eta_k \eta_{k-1}s_0\dots s_k}
\ar[d]_-{ \eta_{k-1}s_0\dots s_k}&
s_ks_{k-1}s_0\dots s_k\ar[r]^-{\lambda_{k-1,k}s_0\dots s_k}
\ar[rd]_-{\lambda_{k-1,k}s_0\dots s_k}^-{\eqref{eq:wdlid}}\ar[ddd]&
s_{k-1}s_ks_0\dots s_k\ar[d]^-{\overline \lambda_{k-1,k}s_0\dots s_k}\\
s_{k-1}s_0\dots s_k\ar[dd]&
\ar@{}[rdd]|-{\mathrm{(YB)}}&
s_{k-1}s_ks_0\dots s_k\ar[d]\\
&&
s_0\dots s_{k-2}(s_{k-1}s_k)^2\ar[d]\\
s_0\dots s_{k-2}s_{k-1}^2s_k\ar[r]^-{\eta_k s_0\dots s_{k-2}s_{k-1}^2s_k}
\ar[d]_-{s_0\dots s_{k-2}\mu_{k-1}s_k}&
s_ks_0\dots s_{k-2}s_{k-1}^2s_k\ar[r]\ar[d]_-{s_ks_0\dots s_{k-2}\mu_{k-1}s_k}
\ar@{}[rd]|-{\hspace{-.8cm}\eqref{eq:wdl}}&
s_0\dots s_{k-2}s_{k-1}^2s_k^2\ar[d]^-{s_0\dots s_{k-2}\mu_{k-1}s_k^2}\\
s_0\dots s_k\ar[r]_-{\eta_k s_0\dots s_k}&
s_ks_0\dots s_k\ar[r]&
s_0\dots s_{k-1}s_k^2
}
$$
we obtain
$\overleftarrow\lambda_{0,\dots,k-2,(k-1,k)}=
\overleftarrow\lambda_{0,\dots,k}.\overleftarrow\lambda_{0,\dots,k-1}s_k=
s_0\dots s_{k-2}\overline\lambda_{k-,k}.
\overleftarrow\lambda_{0,\dots,k}.\overleftarrow\lambda_{0,\dots,k-1}s_k$. 
Combining these identities we conclude the claim.
\end{proof}

\begin{lemma}\label{lem:lambda_n}
In terms of the 2-functors in Theorem \ref{thm:n_w-wreath}, for any
integer $n>1$, the composite 
$$
\xymatrix{
\Wdl^{(n)}(\ocK)\ar[r]^-{C_1}&
\Wdl^{(n-1)}(\ocK)\ar[r]^-{C_1}&
\dots \ar[r]^-{C_1}&
\Wdl(\ocK)}
$$
takes a 0-cell $\{\lambda_{i,j}:s_js_i\to s_i s_j\}_{0\leq i<j\leq n}$ to the
weak distributive law
\begin{equation}\label{eq:lambda_n}
\lambda_{0\dots n-1,n}:=\big(
\xymatrix{
s_n(s_0\dots s_{n-1})\ar[rr]^-{s_n\overline\lambda_{0\dots n-1}}&&
s_ns_0\dots s_{n-1} \ar[r]&
(s_0\dots s_{n-1})s_n
}\big),
\end{equation}
where the unlabelled arrow denotes the unique combination of $\lambda_{i,j}$s
with the given domain and codomain. 
\end{lemma}

\begin{proof}
We proceed by induction in $n$. For $n=2$ the claim follows by Theorem
\ref{thm:n_w-wreath}. Assume that it holds for some $n\geq 2$. Then the
following diagram commutes.   
$$
\xymatrix @C=40pt {
s_0\dots s_n \ar[r]^-{\eta_n s_0\dots s_n}
\ar[d]_-{\overline\lambda_{0\dots n-1}s_n}&
s_ns_0\dots s_n \ar[d]^-{s_n \overline\lambda_{0\dots n-1}s_n}
\ar@/^2pc/[rd]^(.7){\lambda_{0\dots n-1,n} s_n}&
\ar@{}[ld]|(.35){(\ast)}\\
s_0\dots s_n \ar[r]^-{\eta_n s_0\dots s_n}
\ar@/_2.2pc/[rrd]_(.2){\overleftarrow \lambda_{0,\dots, n}}&
s_ns_0\dots s_n \ar[r]&
s_0\dots s_{n-1}s_n ^2\ar[d]^-{s_0\dots s_{n-1}\mu_n}\\
&&s_0\dots s_n}
$$
The region marked by $(\ast)$ commutes by the induction hypothesis. The left
bottom path is equal to $\overline \lambda_{0\dots n}$, see 
\eqref{eq:n_idemp}. Hence we conclude that $\overline \lambda_{0\dots n}$ is
equal to $\overline \lambda_{0\dots n-1,n}$; i.e. the idempotent associated to
the weak distributive law $\lambda_{0\dots n-1,n}$. Using this observation and
\eqref{eq:Ck0} (for the monads $s_0\dots s_{n-1}$, $s_n$ and $s_{n+1}$), in
the top-right path of the following diagram we recognize $\lambda_{0\dots
n,n+1}$.  
$$
\xymatrix @C=40pt {
s_{n+1}s_0\dots s_n\ar[r]^-{s_{n+1}\overline\lambda_{0\dots n}}
\ar[rd]_-{s_{n+1}\overline\lambda_{0\dots n}}&
s_{n+1}s_0\dots s_n\ar[d]^-{s_{n+1}\overline\lambda_{0\dots n-1}s_n}
\ar@/^2.1pc/[rd]^(.7){\lambda_{0\dots n-1,n+1}s_n}&
\ar@{}[ld]|(.35){(\ast)}\\
&s_{n+1}s_0\dots s_n\ar[r]\ar[rd]&
s_0\dots s_{n-1}s_{n+1}s_n\ar[d]^-{s_0\dots s_{n-1}\lambda_{n,n+1}}\\
&&s_0\dots s_{n+1}}
$$
The region marked by $(\ast)$ commutes by the induction hypothesis and the
triangle on the left commutes by \eqref{eq:lambda_n_idemp}.
\end{proof}

Applying \eqref{eq:rarrowintw} and \eqref{eq:rarrow_bar} to the monads
$s_0\dots s_{n-2}$, $s_{n-1}$ and $s_n$, from \eqref{eq:lambda_n} we obtain
the equal expression   
\begin{equation}\label{eq:lambda_n_alt}
\lambda_{0\dots n-1,n}=\big(
\xymatrix{
s_n(s_0\dots s_{n-1})\ar[r]&
s_0\dots s_n\ar[r]^-{\overline\lambda_{0\dots n}}&
(s_0\dots s_{n-1})s_n}\big).
\end{equation}

\begin{lemma}\label{lem:C1}
In terms of the 2-functors in Theorem \ref{thm:n_w-wreath}, for any
positive integer $n$, consider the composite 
\begin{equation}\label{eq:C_1s}
\xymatrix{
\Wdl^{(n)}(\ocK)\ar[r]^-{C_1}&
\Wdl^{(n-1)}(\ocK)\ar[r]^-{C_1}&
\dots \ar[r]^-{C_1}&
\Wdl(\ocK)\ar[r]^-{C_1}&
\Mnd(\ocK)\ .
}
\end{equation}
\begin{itemize}
\item[{(1)}]
It takes a 0-cell $\{\lambda_{i,j}:s_js_i\to s_i
s_j\}_{0\leq i<j\leq n}$ to the monad $(s_0s_1\dots s_n,\overline
\lambda_{01\dots n})$ in $\ocK$, with multiplication $\mu_{0\dots n}$ equal to 
\begin{equation}\label{eq:n_wr_mp}
\xymatrix @C=35pt{
s_0s_1\dots s_ns_0s_1\dots s_n\ar[r]&
s_0^2s_1^2\dots s_n^2\ar[r]^-{\mu_0\mu_1\dots\mu_n}&
s_0s_1\dots s_n \ar[r]^-{\overline \lambda_{01\dots n}}&
s_0s_1\dots s_n 
}
\end{equation}
and unit $\eta_{0\dots n}$ equal to 
\begin{equation}\label{eq:n_wr_unit}
\xymatrix @C=40pt{
1\ar[r]^-{\eta_n\eta_{n-1}\dots\eta_0}&
s_ns_{n-1}\dots s_0\ar[r]&
s_0s_1\dots s_n \ar[r]^-{\overline \lambda_{01\dots n}}&
s_0s_1\dots s_n ,
}
\end{equation}
where the unlabelled arrows denote the unique (by the Yang-Baxter condition)
combinations of $\lambda_{i,j}$s with the given domain and codomain.
\item[{(2)}] It takes a 1-cell $\{\xi_i:s'_iv\to vs_i\}_{0\leq i\leq n}$ to 
$$
\scalebox{.9}{
\xymatrix @C=35pt{
s_0's'_1\dots s'_nv\ar[r]^-{s_0's'_1\dots s'_{n-1}\xi_n}&
s_0's'_1\dots s'_{n-1}vs_n\ar[rr]^-{s_0's'_1\dots s'_{n-2}\xi_{n-1}s_n}&&
\dots\ar[r]^-{\xi_0s_1s_2\dots s_n}&
vs_0s_1\dots s_n\ar[r]^-{v\overline\lambda_{01\dots n}}&
vs_0s_1\dots s_n,
}}
$$
to be denoted by $\xi_{01\dots n}$.
\item[{(3)}]
On the 2-cells it acts as the identity map.
\end{itemize}
\end{lemma}

\begin{proof}
(1) We proceed by induction in $n$. For $n=1$ the claim holds by Paragraph
\ref{as:w_wreath}. Assume that it holds for some $n\geq 1$. Then $\mu_{0\dots 
n+1}$ occurs in the top-right path of the following diagram.
$$
\scalebox{.9}{
\xymatrix{
(s_0\dots s_{n+1})^2\ar[rrrr]^-{s_0\dots s_n\lambda_{0\dots n,n+1}s_{n+1}}
\ar@{=}[d]\ar@{}[rrrrd]|-{\eqref{eq:lambda_n_alt}}&&&&
(s_0\dots s_n)^2s_{n+1}^2\ar@{=}[d]\\
(s_0\dots s_{n+1})^2\ar[r]\ar[rd]&
(s_0\dots s_n)^2s_{n+1}^2\ar[d]\ar@{=}[rr] \ar@{}[rrd]|-{(\ast)}&&
(s_0\dots s_n)^2s_{n+1}^2
\ar[r]^-{\raisebox{10pt}{${}_{s_0\dots s_n\overline\lambda_{0\dots n+1}s_{n+1}}$}}
\ar[d]^-{\mu_{0\dots n}\mu_{n+1}}&
(s_0\dots s_n)^2s_{n+1}^2 \ar[d]_-{\mu_{0\dots n}\mu_{n+1}}\\
&s_0^2\dots s_{n+1}^2\ar[r]_-{\mu_0\dots \mu_{n+1}}&
s_0\dots s_{n+1}\ar[r]_-{\overline\lambda_{0\dots n}s_{n+1}}&
s_0\dots s_{n+1}\ar[r]_-{\overline\lambda_{0\dots n+1}}&
s_0\dots s_{n+1}}}
$$
The region marked by $(\ast)$ commutes by the induction hypothesis and the
bottom-right region commutes by the bilinearity of $\overline\lambda_{0\dots
n+1}= \overline\lambda_{0\dots n,n+1}$, cf. \eqref{eq:wdlid}.
The composite of the last two arrows in the bottom row is equal to
$\overline\lambda_{0\dots n+1}$ by the explicit form in \eqref{eq:n_idemp}. 

Similarly, $\eta_{0\dots n+1}$ occurs in the top-right path of the following
diagram. 
$$
\xymatrix @C=60pt{
A\ar[rrr]^-{\eta_{n+1} \eta_{0\dots n}}\ar@{=}[d]
\ar@{}[rrrd]|-{(\ast)}&&&
s_{n+1}s_0\dots s_n\ar@{=}[d]\\
A\ar[r]^-{\eta_{n+1}\dots \eta_0}&
s_{n+1}\dots s_0\ar[r]\ar[rd]&
s_{n+1}s_0\dots s_n\ar[r]^-{s_{n+1}\overline\lambda_{0\dots n}}
\ar[rd]^-{\lambda_{0\dots n,n+1}}\ar[d]&
s_{n+1}s_0\dots s_n\ar[d]^-{\lambda_{0\dots n,n+1}}\\
&&s_0\dots s_{n+1}\ar[r]_-{\overline\lambda_{0\dots n+1}}
\ar@{}[ru]|(.3){\eqref{eq:lambda_n_alt}}&
s_0\dots s_{n+1}}
$$
The region marked by $(\ast)$ commutes by the induction hypothesis.

Part (2) is easily proved by induction in $n$ and part (3) is trivial.
\end{proof}

\begin{theorem}\label{thm:wreath_assoc}
For any 2-category $\cK$, and any positive integer $n$, the 2-functors
in Theorem \ref{thm:n_w-wreath} give rise to a unique composite 
$$
\xymatrix {
\Wdl^{(n)}(\ocK)\ar[r]^-{C_{k_n}}&
\Wdl^{(n-1)}(\ocK)\ar[r]^-{C_{k_{n-1}}}&
\dots \ar[r]^-{C_{k_2}}&
\Wdl(\ocK)
\ar[r]^-{C_{k_1}}&
\Mnd(\ocK)
}
$$
which is independent of the choice of the index set $k_i\in \{1,2,\dots,i\}$,
for $1\leq i\leq n$.
\end{theorem}

\begin{proof}
We proceed by induction in $n$. For $n=1$ the claim is trivial: there is only
one 2-functor $C_1:\Wdl^{(1)}(\ocK)\equiv \Wdl(\ocK)\to \Wdl^{(0)}(\ocK)\equiv
\Mnd(\ocK)$, that recalled in Paragraph \ref{as:w_wreath}. 
Assume now that the claim holds for some positive integer $n$; i.e. 
$C_1C_{k_2}\dots C_{k_n}$ does not depend on $k_i\in \{1,2,\dots,i\}$, for 
$1\leq i\leq n$. 
With the explicit form of the 2-functor \eqref{eq:C_1s} in Lemma \ref{lem:C1},
and the explicit form of 
$C_k:\Wdl^{(n+1)}(\ocK)\to \Wdl^{(n)}(\ocK)$ in \eqref{eq:Ck0} and
\eqref{eq:Ck1} at hand, the equality $C_kC_1C_1\dots C_1=C_1C_1C_1\dots C_1$
of the $n+1$-fold composites follows by Lemma \ref{lem:n_idemp}, for all
$1\leq k \leq n+1$.  
\end{proof}
 
\begin{aussage}{\bf If idempotent 2-cells split.}
Let us take a 2-category $\cK$ which is locally idempotent complete;
i.e. $\cK\simeq \ocK$. Then we may consider the pseudofunctors 
$$
\xymatrix{
\Wdl^{(n)}(\cK)\ar[r]^-\simeq &
\Wdl^{(n)}(\ocK)\ar[r]^-{C_k}&
\Wdl^{(n-1)}(\ocK)\ar[r]^-\simeq&
\Wdl^{(n-1)}(\cK),
}
$$
for all values of $1\leq k\leq n$. (Choosing the biequivalence $\ocK\to \cK$
by adopting the convention that we split identity 2-cells trivially, i.e. via 
identity 2-cells; they become in fact 2-functors.) Their $n$-fold
iteration is pseudonaturally equivalent to  
$$
\xymatrix{
\Wdl^{(n)}(\cK)\ar[r]^-\simeq &
\Wdl^{(n)}(\ocK)\ar[r]&
\Mnd(\ocK)\ar[r]^-\simeq &
\Mnd(\cK)}
$$
where the unlabelled arrow stands for {\em the} 2-functor in Theorem
\ref{thm:wreath_assoc}. 
This pseudofunctor (or in fact 2-functor with an appropriate choice) takes an
object $\{\lambda_ {i,j}:s_js_i\to s_is_j\}_{0\leq i<j\leq n}$ of
$\Wdl^{(n)}(\cK)$, considered as an object of $\Wdl^{(n)}(\ocK)$, to the image
of the idempotent $\overline\lambda_{0\dots n}$ in \eqref{eq:n_idemp}. This is
regarded as the weak wreath product of the monads $s_0,s_1,\dots, s_n$ in
$\cK$. It is unique -- i.e. the weak wreath product is associative -- up-to an
isomorphism arising from the chosen splittings of the occurring idempotents.
\end{aussage}

\section{Examples from Ising type spin chains}\label{sec:spin}

In this section $\cK:=\mathsf{Vec}$ will be the one-object 2-category (in fact 
bicategory); i.e. monoidal category of vector spaces over a given field
$F$. Thus there is only one 0-cell $\ast$; the 1-cells are the $F$-vector
spaces and the 2-cells are the linear maps. The horizontal composition
(i.e. monoidal product) is given by the tensor product $\ox$ and the vertical
composition is given by the composition of linear maps. Monads are just the
$F$-algebras. We shall make use of the fact that the monoidal category of
vector spaces is {\em symmetric}; the symmetry natural isomorphism (i.e. the
flip map) will be denoted by $\sigma$. Clearly, $\mathsf{Vec}$ is idempotent
complete. 

Our aim is to present an object of $\Wdl^{(n)}(\mathsf{Vec})$ (for any
positive integer $n$) in terms of a finite dimensional weak bialgebra.
We start with recalling the notion of {\em weak bialgebra} from
\cite{Nill}, \cite{BNSz}. 

\begin{definition}
A {\em weak bialgebra} is a vector space $H$ equipped with an algebra
(i.e. monad) structure $\mu:H\ox H \to H$, $\eta:F \to H$ and a coalgebra
(i.e. comonad) structure $\Delta:H \to H\ox H$, $\varepsilon:H \to F$ such
that the following diagrams commute.
$$
\xymatrix @C=20pt @R=20pt {
H^{\ox 2} \ar[r]^-{\Delta\ox \Delta}\ar[d]_-{\mu}&
H^{\ox 4} \ar[rr]^-{H\ox \sigma\ox H}&&
H^{\ox4}  \ar[d]^-{\mu\ox \mu}\\
H\ar[rrr]_-{\Delta}&&&
H^{\ox 2} 
}
$$
$$
\xymatrix @C=35pt @R=15pt{
F\ar[r]^-{\eta\ox \eta}\ar[rd]^-{\eta}\ar[d]_-{\eta\ox \eta}&
H^{\ox 2} \ar[r]^-{\Delta\ox \Delta}&
H^{\ox 4}  \ar[dd]^(.3){H\ox \mu\ox H}&
H^{\ox 3}\ar[r]^-{H\ox \Delta\ox H}\ar[rd]^-{\mu^2}
\ar[dd]_(.7){H\ox \Delta\ox H}&
H^{\ox 4}\ar[r]^-{H\ox \sigma\ox H}&
H^{\ox 4}\ar[d]^-{\mu\ox \mu}\\
H^{\ox 2} \ar[d]_-{\Delta\ox \Delta}&
H\ar[rd]^-{\Delta^2}&&&
H\ar[rd]^-{\varepsilon}&
H^{\ox 2}\ar[d]^-{\varepsilon\ox \varepsilon}\\
H^{\ox 4} \ar[r]_-{H\ox \sigma\ox H}&
H^{\ox 4}\ar[r]_-{H\ox \mu\ox H}&
H^{\ox 3}&
H^{\ox 4}\ar[r]_-{\mu\ox \mu}&
H^{\ox 2}\ar[r]_-{\varepsilon\ox \varepsilon}&
F
}
$$
\end{definition}
This definition can be generalized from $\mathsf{Vec}$ to any braided monoidal
category, see \cite{AAatAl}, \cite{PaSt}. Note that the axioms of a weak
bialgebra are self-dual in the sense that they are closed under the reversing
of the arrows in the representing diagrams.  

\begin{aussage}{\bf Duals of weak bialgebras.}
Whenever the 1-cell underlying a monad in a 2-category possesses a (left or
right) adjoint, this adjoint comes equipped with the canonical structure
of a comonad. Conversely, the adjoint of a comonad is a monad. In the monoidal
category of vector spaces, a 1-cell $H$ -- i.e. a vector space -- possesses a
(left and right) adjoint if and only if it is finite dimensional over $F$; in
which case the adjoint is the linear dual $\hat H:=\mathsf{Hom}(H,F)$; that
is, the vector space of linear maps from $H$ to $F$. In particular, the dual
of a finite dimensional weak bialgebra is both an algebra -- with
multiplication $\hat \mu:=\mathsf{Hom}(\Delta,F):\hat H\ox \hat H \cong 
\mathsf{Hom}(H\ox H,F)\to \hat H$ and unit $\hat \eta:=
\mathsf{Hom}(\varepsilon,F): F\to \hat H$ -- and a coalgebra -- with
comultiplication $\hat \Delta:=\mathsf{Hom}(\mu,F)$ and counit $\hat
\varepsilon:= \mathsf{Hom}(\eta,F)$. That is to say, the (co)algebra structure
of $\hat H$ is defined by the following commutative diagrams
\begin{equation}\label{eq:hat_H}
\scalebox{.8}{
\xymatrix{
\hat H \ox H^{\ox 2}\ar[r]^-{\hat \Delta \ox H^{\ox 2}}
\ar[d]^-{\hat H \ox \sigma}&
\hat H^{\ox 2}\ox H^{\ox 2}\ar[d]_-{\hat H \ox \mathsf{ev}\ox H}
&
\hat H \ar[r]^-{\hat H \ox \eta}\ar[rdd]_-{\hat \varepsilon}&
\hat H \ox H \ar[dd]_-{\mathsf{ev}}
&
\hat H^{\ox 2}\ox H\ar[r]^-{\sigma\ox H}\ar[d]^-{\hat H^{\ox 2}\ox\Delta}&
\hat H^{\ox 2}\ox H\ar[d]_-{\hat \mu \ox H}
&
H\ar[r]^-{\hat \eta \ox H}\ar[rdd]_-{\varepsilon}&
\hat H \ox H\ar[dd]_-{\mathsf{ev}}\\
\hat H \ox H^{\ox 2}\ar[d]^-{\hat H \ox \mu}&
\hat H \ox H \ar[d]_-{\mathsf{ev}}
&
&&
\hat H^{\ox 2}\ox H^{\ox 2}\ar[d]^-{\hat H \ox \mathsf{ev}\ox H}&
\hat H\ox H \ar[d]_-{\mathsf{ev}}\\
\hat H \ox H \ar[r]_-{\mathsf{ev}}&
F
&
&F
&
\hat H \ox H \ar[r]_-{\mathsf{ev}}&
F
&
&F}}
\end{equation}
where $\mathsf{ev}:\hat H \ox H\to F$ stands for the evaluation map (i.e. the
counit of the adjunction $\hat H \dashv H$).
What is more, by self-duality of the weak bialgebra axioms, $\hat H$ is a weak
bialgebra again with the above algebra and coalgebra structures.  
\end{aussage}

\begin{aussage} {\bf The iterated weak wreath product of a finite weak
    bialgebra and its dual.}
In terms of a finite dimensional weak bialgebra $H$, an object of
$\Wdl^{(n)}(\mathsf{Vec})$ is given as follows. If $0\leq i\leq n$ is even,
then let $s_i$ be the algebra underlying $H$ and if $i$ is odd then let $s_i$
be the algebra underlying $\hat H$. If $j-i>1$ then let $\lambda_{i,j}$ be
given by the flip map $\sigma$. If $i$ is odd, then let $\lambda_{i,i+1}$ be
equal to $\lambda$ defined as 
$$
\xymatrix @C=12pt{
H\ox \hat H\ar[r]^-\sigma&
\hat H \ox H \ar[rr]^-{\hat \Delta\ox \Delta}&&
\hat H \ox \hat H \ox H\ox H\ar[rr]^-{\hat H \ox \mathsf{ev} \ox H}&&
\hat H \ox H,
}
$$
and if $i$ is even then let $\lambda_{i,i+1}$ be equal to $\hat \lambda$ given
by 
$$
\xymatrix @C=12pt{
\hat H\ox H\ar[r]^-\sigma&
H \ox \hat H \ar[rr]^-{\Delta\ox \hat \Delta}&&
H \ox H \ox \hat H\ox \hat H\ar[rr]^-{H \ox \sigma \ox \hat H}&&
H \ox \hat H \ox H \ox \hat H \ar[rr]^-{H \ox \mathsf{ev} \ox \hat H}&&
H \ox \hat H.
}
$$
The symmetry $\sigma:X \ox Y \to Y \ox X$, for $X,Y\in \{H,\hat H\}$, is a
distributive law hence a weak distributive law. We show that $\lambda:H\ox
\hat H \to \hat H \ox H$ is a weak distributive law. The morphism 
$$
\xi:=\big(
\xymatrix{
\hat H \ox H\ar[r]^-{\hat \Delta \ox H}&
\hat H^{\ox 2}\ox H \ar[r]^-{\hat H \ox \mathsf{ev}}&
\hat H}\big)
$$
is an associative (and evidently unital) action in the sense of commutativity
of
\begin{equation}\label{eq:lambda_associ}
\xymatrix @C=45pt @R=18pt {
\hat H \ox H^{\ox 2}\ar[rr]^-{\hat \Delta \ox H^{\ox 2}}
\ar[dd]_-{\hat H \ox \sigma}\ar[rd]^-{\hat\Delta\ox H^{\ox 2}}&&
\hat H^{\ox 2}\ox H^{\ox 2}\ar[r]^-{\hat H \ox \mathsf{ev}\ox H}
\ar[d]^-{\hat \Delta \ox \hat H \ox H^{\ox 2}}&
\hat H \ox H\ar[d]^-{\hat \Delta \ox H}\\
&\hat H^{\ox 2}\ox H^{\ox 2}\ar[r]^-{\hat H \ox \hat \Delta\ox H^{\ox 2}}
\ar[d]^-{\hat H^{\ox 2}\ox \sigma}&
\hat H^{\ox 3}\ox H^{\ox 2}\ar[r]^-{\hat H^{\ox 2}\ox \mathsf{ev}\ox H}&
\hat H^{\ox 2}\ox H \ar[dd]^-{\hat H \ox \mathsf{ev}}\\
\hat H \ox H^{\ox 2}\ar[r]^-{\hat \Delta \ox H^{\ox 2}}
\ar[d]_-{\hat H \ox \mu}&
\hat H^{\ox 2}\ox H^{\ox 2}\ar[d]^-{\hat H^{\ox 2}\ox \mu}\\
\hat H \ox H \ar[r]_-{\hat \Delta \ox H}&
\hat H^{\ox 2}\ox H \ar[rr]_-{\hat H \ox \mathsf{ev}}&&
\hat H}
\end{equation}
where the bottom-right region commutes by the first identity in
\eqref{eq:hat_H}. In terms of $\xi$,
$$
\lambda =\big(
\xymatrix{
H\ox \hat H\ar[r]^-{\sigma}&
\hat H \ox H \ar[r]^-{\hat H \ox \Delta}&
\hat H \ox H^{\ox 2}\ar[r]^-{\xi\ox H}&
\hat H \ox H
}\big).
$$
Using this form of $\lambda$, its compatibility with the multiplication of $H$
follows by commutativity of the diagram below.
$$
\scalebox{.85}{\xymatrix{
H^{\ox 2}\ox \hat H\ar[rr]^-{H\ox \sigma}\ar[rrd]^-{H\ox \Delta \ox \hat H}
\ar[rddd]^-{\sigma_{H\ox H,\hat H}}\ar[dddd]^-{\mu \ox \hat H}&&
H\ox \hat H \ox H\ar[r]^-{H\ox \hat H \ox \Delta}&
H\ox \hat H \ox H^{\ox 2}\ar[rr]^-{H \ox \xi \ox H}
\ar[d]^-{\sigma_{H,\hat H\ox H}\ox H}&&
H\ox \hat H \ox H\ar[d]_-{\sigma\ox H}\\
&& 
H^{\ox 3}\ox \hat H\ar[ru]_-{\ H\ox \sigma_{H\ox H,\hat H}}
\ar[d]_-{\Delta \ox H^{\ox 2}\ox \hat H}&
\hat H \ox H^{\ox 3}\ar[rr]^-{\xi \ox H^{\ox 2}}
\ar[d]^-{\hat H \ox H \ox \Delta \ox H}&&
\hat H \ox H^{\ox 2}\ar[d]_-{\hat H \ox \Delta \ox H}\\
&& 
H^{\ox 4}\ox \hat H\ar[d]_-{\sigma_{H^{\ox 4},\hat H}}&
\hat H \ox H^{\ox 4}\ar[rr]^-{\xi \ox H^{\ox 3}}
\ar[d]^-{\hat H \ox \sigma\ox H^{\ox 2}}
\ar@{}[rrd]|-{\eqref{eq:lambda_associ}}&&
\hat H \ox H^{\ox 3}\ar[d]_-{\xi\ox H^{\ox 2}}\\
&
\hat H \ox H^{\ox 2}\ar[r]^-{\hat H \ox \Delta \ox \Delta}
\ar[d]^-{\hat H \ox \mu}&
\hat H \ox H^{\ox 4}\ar[r]^-{\hat H \ox H \ox \sigma\ox H}
\ar[ru]^-{\hat H \ox \sigma_{H\ox H,H} \ox H\ }&
\hat H \ox H^{\ox 4}\ar[r]^-{\hat H \ox \mu \ox H^{\ox 2}}&
\hat H \ox H^{\ox 3}\ar[r]^-{\xi \ox H^{\ox 2}}\ar[d]^-{\hat H \ox H \ox \mu}&
\hat H \ox H^{\ox 2}\ar[d]_-{\hat H \ox \mu}\\
H \ox \hat H \ar[r]_-{\sigma}&
\hat H \ox H\ar[rrr]_-{\hat H \ox \Delta}&&&
\hat H \ox H^{\ox 2}\ar[r]_-{\xi \ox H}&
\hat H \ox H}}
$$
The region at the middle of the bottom row commutes by the first weak
bialgebra axiom. Symmetrically, in terms of $\zeta:=(\mathsf{ev}\ox H).(\hat H
\ox \Delta)$,
we can write $\lambda=(\hat H \ox \zeta).(\hat \Delta \ox H).\sigma$. With
this form of $\lambda$ at hand, its compatibility with the multiplication of
$\hat H$ follows symmetrically. It remains to check the weak unitality
condition \eqref{eq:wunit}. For that consider the (idempotent) morphism
$$
\bar\varepsilon_s:=\big(\xymatrix{
H\ar[r]^-{\eta \ox H}&
H^{\ox 2}\ar[r]^-{\Delta\ox H}&
H^{\ox 3}\ar[r]^-{H\ox \mu}&
H^{\ox 2}\ar[r]^-{H\ox \varepsilon}&
H}\big).
$$
Recall from \cite{BoCaJa:wba} (equations (4) and (8), respectively) that the
following diagrams involving $\bar\varepsilon_s$ commute.
\begin{equation}\label{eq:BCJ}
\xymatrix{
H^{\ox 2}\ar[r]^-{\Delta \ox H}\ar[d]_-{H\ox \bar\varepsilon_s}&
H^{\ox 3}\ar[r]^-{H\ox \mu}&
H^{\ox 2}\ar[d]^-{H\ox \varepsilon}
&&
H \ar[r]^-{\eta\ox H}\ar[d]_-\Delta&
H^{\ox 2}\ar[r]^-{\Delta \ox H}&
H^{\ox 3}\ar[d]^-{H\ox \mu}\\
H^{\ox 2} \ar[rr]_-\mu&&
H
&&
H^{\ox 2}\ar[rr]_-{\bar\varepsilon_s\ox H}&&
H^{\ox 2}}
\end{equation}
Using the definitions in \eqref{eq:hat_H}, the first identity in
\eqref{eq:BCJ} is equivalent to 
$\xi.(\hat H \ox \bar \varepsilon_s)=
\hat \mu.(\hat H \ox \xi).(\hat H \ox \hat \eta \ox H)$, implying
commutativity of the bottom-right region in 
$$
\xymatrix{
\hat H \ox H \ar[r]^-{\eta \ox \hat H \ox H}\ar@{=}[d] 
\ar@/_1pc/[rr]_-{\hat H \ox \eta\ox H}&
H \ox \hat H \ox H \ar[r]^-{\sigma\ox H}&
\hat H \ox H^{\ox 2}\ar[r]^-{\hat H \ox \Delta \ox H}
\ar@{}[d]|-{\eqref{eq:BCJ}}&
\hat H \ox H^{\ox 3}\ar[r]^-{\xi \ox H^{\ox 2}}
\ar[d]^-{\hat H \ox H \ox \mu}&
\hat H \ox H^{\ox 2}\ar[d]^-{\hat H \ox \mu}\\
\hat H \ox H\ar[rr]^-{\hat H \ox \Delta}\ar[d]_-{\hat H \ox H\ox \hat \eta}
\ar[rd]^-{\hat H \ox \hat \eta\ox H}&&
\hat H \ox H^{\ox 2}\ar[r]^-{\hat H \ox \bar \varepsilon_s\ox H}
\ar[d]^-{\hat H \ox \hat \eta\ox H^{\ox 2}}
&
\hat H \ox H^{\ox 2}\ar[r]^-{\xi \ox H}&
\hat H \ox H \ar@{=}[d]\\
\hat H \ox H \ox \hat H \ar[r]_-{\hat H \ox \sigma}&
\hat H^{\ox 2}\ox H\ar[r]_-{\hat H^{\ox 2}\ox \Delta}&
\hat H^{\ox 2}\ox H^{\ox 2}\ar[r]_-{\hat H \ox \xi \ox H}&
\hat H^{\ox 2}\ox H\ar[r]_-{\hat \mu \ox H}&
\hat H \ox H\ .}
$$
Any path in this diagram yields an alternative expression of
the idempotent $\overline \lambda:\hat H \ox H \to \hat H \ox H$, proving
that $\lambda$ is a weak distributive law. By symmetrical considerations so is
$\hat \lambda$.  

With some routine computations using the weak bialgebra axioms, one checks
that $\overline \lambda$ is equal to the identity map -- i.e. $\lambda$ is a
distributive law in the strict sense -- if and only if $\Delta.\eta=\eta\ox
\eta$; i.e. $H$ is a bialgebra in the strict sense. 

Our next task is to check the Yang-Baxter conditions. The symmetry operators
among themselves obey the Yang-Baxter condition, hence for $\{i,j,k\}$ such
that $j-i>1 $ and $k-j> 1$ we are done. For $\{i-1,i,j\}$ and $\{i,j,j+1\}$,
such that $j-i> 1$, the Yang-Baxter conditions follow by naturality of the
symmetry. So we are left with the case $\{i-1,i,i+1\}$. Assume first that $i$
is odd. Then the Yang-Baxter condition follows by commutativity of
$$
\xymatrix @R=15pt {
H\ox \hat H \ox H \ar[r]^-{\sigma\ox H}\ar[dd]^-{H\ox \sigma}&
\hat H \ox H^{\ox 2}\ar[r]^-{\hat \Delta \ox \Delta \ox H}
\ar[dd]^-{\sigma_{\hat H \ox H,H}}&
\hat H^{\ox 2} \ox H^{\ox 3}\ar[r]^-{\hat H \ox \mathsf{ev}\ox H^{\ox 2}}
\ar[dd]^-{\sigma_{\hat H^{\ox 2} \ox H^{\ox 2},H}}&
\hat H \ox H^{\ox 2} \ar[d]_-{\hat H \ox \sigma}\\
&&&
\hat H \ox H^{\ox 2} \ar[d]_-{\sigma\ox H}\\
H^{\ox 2}\ox \hat H \ar[r]^-{\sigma_{H,H\ox \hat H}}
\ar[d]^-{H\ox \Delta \ox \hat \Delta}&
H \ox \hat H \ox H 
\ar[r]^-{H\ox \hat \Delta \ox \Delta}
\ar[d]^(.4){\Delta \ox \hat \Delta \ox H}&
H\ox \hat H^{\ox 2}\ox H^{\ox 2}\ar[r]^-{\ H\ox \hat H \ox \mathsf{ev}\ox H}
\ar[d]^-{\Delta \ox \hat \Delta \ox \hat H \ox H^{\ox 2}}&
H\ox \hat H \ox H \ar[d]_-{\Delta \ox \hat \Delta \ox H}\\
H^{\ox 3}\ox \hat H^{\ox 2}
\ar[d]^-{H^{\ox 2}\ox \mathsf{ev}.\sigma \ox \hat H}&
H^{\ox 2}\ox \hat H^{\ox 2}\ox H
\ar[r]^-{\raisebox{5pt}{${}_{H^{\ox 2}\ox \hat H \ox \hat \Delta \ox \Delta}$}}
\ar[dd]^-{H\ox \mathsf{ev}.\sigma \ox \hat H \ox H}&
H^{\ox 2}\ox \hat H^{\ox 3}\ox H^{\ox 2}
\ar[dd]^-{H\ox \ox \mathsf{ev}.\sigma \ox \hat H^{\ox 2} \ox H^{\ox 2}}&
H^{\ox 2}\ox \hat H^{\ox 2}\ox H
\ar[dd]_(.7){H\ox \mathsf{ev}.\sigma \ox \hat H \ox H}\\
H^{\ox 2} \ox \hat H \ar[rd]^-{\sigma_{H,H\ox \hat H}}
\ar[d]^-{\sigma \ox \hat H}\\
H^{\ox 2}\ox \hat H\ar[r]_-{H\ox \sigma}&
H\ox \hat H \ox H \ar[r]_-{H\ox \hat \Delta \ox \Delta}&
H\ox \hat H^{\ox 2} \ox H^{\ox 2}\ar[r]_-{H \ox \hat H \ox \mathsf{ev}\ox H\ }&
H\ox \hat H \ox H\ .}
$$
The case when $i$ is even is treated symmetrically. This proves that the
construction in this paragraph yields an object in
$\Wdl^{(n)}(\mathsf{Vec})$ (which is an object of $\Mnd^{n+1}(\mathsf{Vec})$
if and only if $H$ is a bialgebra in the strict sense). Hence by Theorem
\ref{thm:wreath_assoc} there is a corresponding weak wreath product monad
(i.e. $F$-algebra) given as the image of the idempotent \eqref{eq:n_idemp}. 
Since $\sigma$ is unital, one obtains the following explicit forms of this
idempotent. If $n$ is odd, then it comes out as
$$
\xymatrix @C=40pt {
(H\ox \hat H)^{\ox \frac{n+1}2}
\ar[r]^-{\overline{\hat \lambda}^{\ \ox \frac{n+1}2}}&
H\ox (\hat H \ox H)^{\ox \frac{n-1}2}\ox \hat H 
\ar[rr]^-{H\ox \overline \lambda^{\ \ox \frac{n-1}2}\ox \hat H}&&
(H\ox \hat H)^{\ox \frac{n+1}2}
}
$$
and if $n$ is even, then 
$$
\xymatrix @C=40pt{
(H\ox \hat H)^{\ox \frac{n}2}\ox H
\ar[r]^-{\overline{\hat \lambda}^{\ \ox \frac{n}2}\ox H}& 
H\ox (\hat H \ox H)^{\ox \frac{n}2}
\ar[r]^-{H\ox \overline \lambda^{\ \ox \frac{n}2}}&
(H\ox \hat H)^{\ox \frac{n}2}\ox H.
}
$$
If $H$ is a bialgebra in the strict sense (e.g. it is the linear span of a
finite group), then these idempotents become identity maps and so the above
weak wreath products reduce wreath products in the strict sense. 
\end{aussage}

In the quantum spin chains in \cite{NilSzlWie:Jones} and \cite{Bohm:PhD},
where the spins take their values in a dual pair of finite dimensional weak
Hopf algebras, this (n+1)-ary weak wreath product is regarded as the algebra
of  observable quantities localized in the interval $[0,n]$ of the one
dimensional lattice. In particular, in spin chains built on dual pairs of finite
dimensional Hopf algebras (e.g. pairs of a finite group algebra and the
algebra of linear functions on this group), the observable algebra is a proper 
(n+1)-ary wreath product. In the classical Ising model -- where the spins only
have `up' and `down' positions -- these dual Hopf algebras are both isomorphic
to the linear span of the sign group $\mathbb Z(2)$.

\section{A fully faithful embedding} \label{sec:embed}

In this section we show that, for any 2-category $\cK$, and any non-negative
integer $n$, $\Wdl^{(n)}(\ocK)$ admits a fully faithful embedding into the
power 2-category $\Mnd(\ocK)^{\two^{n+1}}$. Whenever idempotent 2-cells in
$\cK$ split, this gives rise to a fully faithful embedding $\Wdl^{(n)}(\cK)
\to \Mnd(\cK)^{\two^{n+1}}$. If in addition $\cK$ admits Eilenberg-Moore
objects, this amounts to a fully faithful embedding $\Wdl^{(n)}(\cK)\to
\cK^{\two^{n+1}}$.  

\begin{aussage} \label{as:K^2n}
{\bf The 2-category $\cK^{\two ^n}$.}
The 2-category $\two$ has two 0-cells $0$ and $1$; an only non-identity 1-cell
$1\to 0$; and all of its 2-cells are identities. For any 2-category $\cK$,
there is a 2-category $\cK^\two$ of 2-functors $\two \to \cK$, 2-natural
transformations and modifications. Iteratively, for $n>1$ we define $\cK^{\two^
n}$ as $(\cK^{\two^{n-1}})^\two$. That is, $\cK^{\two^n}$ is isomorphic to the
2-category of 2-functors from the $n$-fold Cartesian product $\two \times
\dots \times \two$ to $\cK$, 2-natural transformations and modifications. An
explicit description is given as follows. The 0-cells are the $n$ dimensional
oriented cubes whose 2-faces are commutative squares of 1-cells in $\cK$. A
1-cell from an $n$-cube of edges 
$\{v_{{\underline p},{\underline q}}: 
A_{\underline p} \to A_{\underline q} \}$ 
to $\{v'_{{\underline p},{\underline q}}: 
A'_{\underline p} \to A'_{\underline q} \}$ 
consists of 1-cells 
$\{u_{\underline p}:A_{\underline p} \to A'_{\underline p}\}$
 in $\cK$ such that the $n+1$-cube 
$\{v_{{\underline p},{\underline q}}:
 A_{\underline p} \to A_{\underline q}, 
u_{\underline p}:A_{\underline p} \to A'_{\underline p}, 
v'_{{\underline p},{\underline q}}: 
A'_{\underline p} \to A'_{\underline q} \}$ 
is commutative. That is, for all values of ${\underline p}$ and
 ${\underline q}$, 
$v'_{{\underline p},{\underline q}}.u_{\underline p}=
u_{\underline q}.v_{{\underline p},{\underline q}}$. Finally, 2-cells consist
of 2-cells $\omega_{\underline p}:u_{\underline p} \to \tilde u_{\underline
  p}$ in $\cK$ such that $ v'_{{\underline p},{\underline q}}.\omega_{\underline p}=
\omega_{\underline q}.v_{{\underline p},{\underline q}}$.
\end{aussage}

In Cartesian coordinates, the vertices of an $n$-cube can be labelled by the
elements $\underline p=(p_1,\dots,p_n)$ of the set $\{0,1\}^{\times n}$. 
Sometimes we represent $\underline p\in \{0,1\}^{\times n}$ by listing those
values of $i$ for which $p_i=1$. For example, $12=(1,1,0,\dots,0)$,
$3=(0,0,1,0,\dots,0)$, etc..  
The $n$-cube has an edge $\underline p \to \underline q$ if and only if there
is some integer $1\leq i \leq n$ such that $p_j=q_j$ for all $j\neq i$,
$p_i=0$ and $q_i=1$. We denote this situation by $\underline q=\underline
p+i$.  
For $\underline p,\underline q\in \{0,1\}^{\times n}$, we say that $\underline
p<\underline q$ if, for any $1\leq i,j\leq n$, the equality $p_iq_j=1$ implies
$i<j$. For $\underline p<\underline q$ we define $\underline p+\underline q\in 
\{0,1\}^{\times n}$ putting $(\underline p+\underline q)_i:=p_i+q_i$. We
denote by $\underline 0:=(0,0,\dots,0)$ and $\underline 1:=(1,1,\dots,1)$ the
constant elements of $\{0,1\}^{\times n}$. 

The construction of the promised 2-functor $\Wdl^{(n)}(\ocK)\to
\Mnd(\ocK)^{\two^{n+1}}$ relies on a few lemmas below. A routine computation
proves the first one: 

\begin{lemma} \label{lem:A}
For any object $\{\lambda_{i,j}:s_js_i\to s_is_j\}_{0\leq i<j\leq 2}$ of
$\Wdl^{(2)}(\ocK)$, there is a homomorphism of monads in $\ocK$,
$$
\varphi_{0,2}^1:= \overline \lambda_{012}. s_0\eta_1s_2:
(s_0s_2,\overline \lambda_{02})\to 
(s_0s_1s_2,\overline\lambda_{012}).
$$
\end{lemma}

This means that there is a 1-cell
$((A,A),\varphi_{0,2}^1):(A,(s_0s_1s_2,\overline\lambda_{012})) \to 
(A,(s_0s_2,\overline \lambda_{02}))$ in $\Mnd(\ocK)$ (where $A$ is the object
underlying the monads $(A,s_i)$).

\begin{lemma}\label{lem:B}
For any object $\{\lambda_{i,j}:s_js_i\to s_is_j\}_{0\leq i<j\leq 4}$ of
$\Wdl^{(4)}(\ocK)$, the morphisms as in Lemma \ref{lem:A} constitute a
commutative diagram in $\ocK$:  
$$
\xymatrix{
(s_0s_2s_4,\overline\lambda_{024})
\ar[r]^-{\varphi_{0,24}^1}\ar[d]_-{\varphi_{02,4}^3}&
(s_0s_1s_2s_4,\overline\lambda_{0124})\ar[d]^-{\varphi_{012,4}^3}\\
(s_0s_2s_3s_4,\overline\lambda_{0234})\ar[r]_-{\varphi_{0,234}^1}&
(s_0s_1s_2s_3s_4,\overline\lambda_{01234})\ .}
$$
\end{lemma}

\begin{proof}
In view of Lemma \ref{lem:A}, both paths around the diagram are equal to
$\overline \lambda_{01234}.$ $s_0\eta_1s_2\eta_3s_4$. 
\end{proof}

\begin{lemma}\label{lem:C}
For any 1-cell $\{\xi_i:s'_i v \to vs_i\}_{0\leq i\leq 2}$ in
$\Wdl^{(2)}(\ocK)$, the morphisms in Lemma \ref{lem:A} induce a commutative
square   
$$
\xymatrix{
(A,(s_0s_1s_2,\overline\lambda_{012}))\ar[rr]^-{((v,\overline v),\xi_{012})}
\ar[d]_-{((A,A),{\varphi}_{0,2}^1)}&&
(A',(s'_0s'_1s'_2,\overline{\lambda'}_{012}))
\ar[d]^-{((A',A'),{\varphi'}_{0,2}^1)}\\
(A,(s_0s_2,\overline \lambda_{02}))\ar[rr]_-{((v,\overline v),\xi_{02})}&&
(A',(s'_0s'_2,\overline {\lambda'}_{02}))}
$$
in $\Mnd(\ocK)$.
\end{lemma}

\begin{proof}
Both paths around the diagram are computed to be equal to
$$
v\overline \lambda_{012}.vs_0\eta_1s_2.\xi_0s_2.s'_0\xi_2=
v\overline \lambda_{012}.\xi_0s_1s_2.s'_0\xi_1s_2.s'_0s'_1\xi_2.
\overline {\lambda'}_{012}v.s'_0\eta'_1s'_2v,
$$
where the last equality follows by the unitality of $\xi_i$ and the
normalization property $\xi_{012}.\overline {\lambda'}_{012}v=\xi_{012}$.
\end{proof}

\begin{aussage} {\bf A 2-functor $\Wdl^{(n)}(\ocK)\to
\Mnd(\ocK)^{\two^{n+1}}$.} \label{as:embed} 
For any any non-negative integer $n$, an any $\underline 0\neq \underline
p=(p_0,\dots,p_n)\in \{0,1\}^{\times(n+1)}$, taking those values of $0\leq
i\leq n$ for that $p_i=1$, defines a 2-functor $\Wdl^{(n)}(\ocK)\to
\Wdl^{(-1+\sum_i p_i)}(\ocK)$. Composing it with the unique iterated weak
wreath product 2-functor $\Wdl^{(-1+\sum_i p_i)}(\ocK)\to \Mnd(\ocK)$ in Section
\ref{sec:iterate}, yields a 2-functor $\Wdl^{(n)}(\ocK)\to \Mnd(\ocK)$. Denote
the image of   
$$
\xymatrix @C=10pt{
\{\lambda_{i,j}:s_js_i\to s_is_j\}_{0\leq i<j\leq n}
\ar@/^2pc/[rr]^-{\{\xi_i:s'_iv\to vs_i\}_{0\leq i \leq n}}
\ar@/_2pc/[rr]_-{\{\xi'_i:s'_iv'\to v's_i\}_{0\leq i \leq n}}&
\Downarrow\omega&
\{\lambda'_{i,j}:s'_js'_i\to s'_is'_j\}_{0\leq i<j\leq n}
}
$$
under it by 
$$
\xymatrix @C=10pt{
(A,(s_{\underline p},\overline\lambda_{\underline p}))
\ar@/^2pc/[rr]^-{\xi_{\underline p}:s'_{\underline p}v\to vs_{\underline p}}
\ar@/_2pc/[rr]_-{\xi'_{\underline p}:s'_{\underline p}v'\to v's_{\underline p}}&
\Downarrow\omega&
(A',(s'_{\underline p},\overline{\lambda'}_{\underline p}))\ .
}
$$
Then $(A,s_{\underline p})$ is the weak wreath product of those demimonads
$(A,s_i)$ for which $p_i=1$. For $\underline 0\in \{0,1\}^{\times(n+1)}$, put
$(A,s_{\underline 0}):=(A,A)$ and $\xi_{\underline 0}:=v$.

By Lemmas \ref{lem:A} and \ref{lem:B}, for any object
$\{\lambda_{i,j}:s_js_i\to s_is_j\}_{0\leq i<j\leq n}$ of $\Wdl^{(n)}(\ocK)$,
for any $\underline p\in \{0,1\}^{\times (n+1)}$, and for any $0\leq i,j\leq
n$ such that $p_i=p_j=0$, there is a commutative square   
$$
\xymatrix @C=40pt{
(A,(s_{\underline p+i+j},\overline\lambda_{\underline p+i+j}))
\ar[r]^-{((A,A),\varphi_{\underline p+j}^i)}
\ar[d]_-{((A,A),\varphi_{\underline p+i}^j)}&
(A,(s_{\underline p+j},\overline\lambda_{\underline p+j}))
\ar[d]^-{((A,A),\varphi_{\underline p}^j)}\\
(A,(s_{\underline p+i},\overline\lambda_{\underline p+i}))
\ar[r]_-{((A,A),\varphi_{\underline p}^i)}&
(A,(s_{\underline p},\overline\lambda_{\underline p}))}
$$
in $\Mnd(\ocK)$.
Such squares constitute a commutative $n+1$-cube in $\Mnd(\ocK)$.

By Lemma \ref{lem:C}, for any 1-cell $\{\xi_i:s'_iv\to vs_i\}_{0\leq i \leq
n}$ of $\Wdl^{(n)}(\ocK)$, there is a commutative square 
$$
\xymatrix{
(A,(s_{\underline p+i},\overline\lambda_{\underline p+i}))
\ar[rr]^-{((v,\overline v),\xi_{\underline p+i})}
\ar[d]_{((A,A),\varphi_{\underline p}^i)}&&
(A',(s'_{\underline p+i},\overline{\lambda'}_{\underline p+i})) 
\ar[d]^{((A',A'),{\varphi'}_{\underline p}^i)}\\
(A,(s_{\underline p},\overline\lambda_{\underline p}))
\ar[rr]_-{((v,\overline v),\xi_{\underline p})}&&
(A',(s'_{\underline p},\overline{\lambda'}_{\underline p}))}
$$
in $\Mnd(\ocK)$. Hence the 1-cells $((v,\overline v),\xi_{\underline p}):
(A,(s_{\underline p},\overline\lambda_{\underline p}))
\to (A',(s'_{\underline p},\overline{\lambda'}_{\underline p}))$ 
constitute a 1-cell in $\Mnd(\ocK)^{\two^{n+1}}$.

Finally, for a 2-cell $\omega:\{\xi_i:s'_iv\to vs_i\}_{0\leq i \leq n} \to
\{\xi'_i:s'_iv'\to v' s_i\}_{0\leq i \leq n}$ in
$\Wdl^{(n)}(\ocK)$, $\omega$ is a 2-cell $((v,\overline v),\xi_{\underline
p})\to ((v',\overline{v'}),\xi'_{\underline p})$ in
$\Mnd(\ocK)$, for any $\underline p\in \{0,1\}^{\times (n+1)}$, which
constitutes evidently a 2-cell in $\Mnd(\ocK)^{\two^{n+1}}$.

The above maps define the stated 2-functor $\Wdl^{(n)}(\ocK)\to
\Mnd(\ocK)^{\two^{n+1}}$. 
\end{aussage}

\begin{theorem} \label{thm:fullyfaithful}
For any 2-category $\cK$, and any non-negative integer $n$, the 2-functor
$\Wdl^{(n)}(\ocK)\to \Mnd(\ocK)^{\two^{n+1}}$ in Paragraph \ref{as:embed} is
fully faithful.  
\end{theorem}

\begin{proof}
Faithfulness is obvious. In order to prove fullness on the 1-cells, take a
1-cell $\{\zeta_{\underline p}:s'_{\underline p}v\to
vs_{\underline p}\}_{{\underline p}\in \{0,1\}^{\times(n+1)}}$ in
$\Mnd(\ocK)^{\two^{n+1}}$ between objects arising from 0-cells
$\{\lambda_{i,j}:s_js_i\to s_is_j\}_{0\leq i<j\leq n}$ and
$\{\lambda'_{i,j}:s'_js'_i\to s'_is'_j\}_{0\leq i<j\leq n}$ of
$\Wdl^{(n)}(\ocK)$. This includes in particular 1-cells
$\xi_i:=\zeta_i:s'_iv\to vs_i$ in $\Mnd(\ocK)$. We 
claim that for $i\in \{0,\dots,n\}$ they constitute a 1-cell in
$\Wdl^{(n)}(\ocK)$ and each $\zeta_{\underline p}$ is equal to their weak
wreath product. 
By commutativity of the squares 
\begin{equation*}
\xymatrix @C=30pt {
(A,(s_is_j,\overline \lambda_{ij})) 
\ar[r]^-{((v,\overline v), \zeta_{ij})}
\ar[d]_-{((A,A),\varphi_i^j)}&
(A',(s'_is'_j,\overline {\lambda'}_{ij}))
\ar[d]^-{((A',A'),{\varphi'}_i^j)}\\
(A,(s_i,\overline s_i))\ar[r]_-{((v,\overline v),\xi_i)}&
(A',(s'_i,\overline s'_i))}\quad
\xymatrix  @C=30pt{
(A,(s_is_j,\overline \lambda_{ij})) 
\ar[r]^-{((v,\overline v), \zeta_{ij})}
\ar[d]_-{((A,A),\varphi_j^i)}&
(A',(s'_is'_j,\overline {\lambda'}_{ij}))
\ar[d]^-{((A',A'),{\varphi'}_j^i)}\\
(A,(s_j,\overline s_j))\ar[r]_-{((v,\overline v),\xi_j)}&
(A',(s'_j,s'_j))}
\end{equation*}
in $\Mnd(\ocK)$, we conclude the commutativity of the diagrams
\begin{equation}\label{eq:xi_ij}
\xymatrix{
s'_iv\ar[d]_-{\xi_i}\ar[r]^-{s'_i\eta'_jv}&
s'_is'_j v \ar[r]^-{\overline{\lambda'}_{ij}v}\ar[rd]^-{ \zeta_{ij}}&
s'_is'_j v\ar[d]^-{ \zeta_{ij}}\\
vs_i\ar[r]_-{vs_i\eta_j}&
vs_is_j\ar[r]_-{v\overline \lambda_{ij}}&
vs_is_j}\qquad
\xymatrix{
s'_jv\ar[d]_-{\xi_j}\ar[r]^-{\eta'_is'_jv}&
s'_is'_j v \ar[r]^-{\overline{\lambda'}_{ij}v}\ar[rd]^-{ \zeta_{ij}}&
s'_is'_j v\ar[d]^-{ \zeta_{ij}}\\
vs_j\ar[r]_-{v\eta_i s_j}&
vs_is_j\ar[r]_-{v\overline \lambda_{ij}}&
vs_is_j}
\end{equation}
in $\cK$. Since $((v,\overline v),\zeta_{ij})$ is a 1-cell in $\Mnd(\ocK)$,
the following diagram commutes
\begin{equation}\label{eq:zeta_wreath}
\xymatrix @R=18pt{
s'_is'_jv\ar[rr]^-{s'_i\xi_j}\ar[ddd]^-{s'_i\eta'_is'_jv}
\ar@/_3pc/[ddddd]_-{\overline {\lambda'}_{ij}v}
\ar@{}[rrd]|-{\eqref{eq:xi_ij}}&&
s'_ivs_j\ar[rr]^-{\xi_is_j}\ar[d]^-{s'_iv\eta_is_j}&&
vs_is_j\ar[d]_-{vs_i\eta_is_j}
\ar@/^3pc/[ddddd]^-{v\overline \lambda_{ij}}\\
&&s'_ivs_is_j\ar[rr]^-{\xi_is_is_j}
\ar[rdd]^-{s'_i\eta'_jvs_is_j}\ar[ldd]_-{s'_iv\overline \lambda_{ij}}
\ar@{}[rrd]|-{\eqref{eq:xi_ij}}&&
vs_i^2s_j\ar[d]_-{vs_i\eta_js_is_j}\\
&&&&v(s_is_j)^2\ar[d]_-{v\overline\lambda_{ij}s_is_j}\\
s_i^{\prime 2}s'_jv\ar[r]^-{s'_i \zeta_{ij}}\ar[d]^(.4){s'_i\eta'_js'_is'_jv}&
s'_ivs_is_j\ar[d]^-{s'_i\eta'_jvs_is_j}&&
s'_is'_jvs_is_j\ar[r]^-{ \zeta_{ij}s_is_j}&
v(s_is_j)^2\ar[dd]_-{v\mu_{ij}}\\
(s'_is'_j)^2v\ar[r]^-{s'_is'_j \zeta_{ij}}\ar[d]^-{\mu'_{ij}v}&
s'_is'_jvs_is_j\ar[rr]^-{ \zeta_{ij}s_is_j}&&
v(s_is_j)^2\ar[rd]^-{v\mu_{ij}}\\
s'_is'_jv\ar[rrrr]_-{ \zeta_{ij}}&&&&
vs_is_j\ ,}
\end{equation}
where we denoted $\mu_{ij}=\mu_i\mu_j.s_i\lambda_{i,j}s_j$.
Since $ \zeta_{ij}.\overline{\lambda'}_{ij}v= \zeta_{ij}$, this
says that $ \zeta_{ij}$ is equal to $v\overline
\lambda_{ij}.\xi_is_j.s'_i\xi_j$; that is, $\zeta_{ij}$ is the weak wreath
product of $\xi_i$ and $\xi_j$. With this expression of $ \zeta_{ij}$ at
hand, also the following diagram commutes.  
$$
\xymatrix @R=18pt{
s'_js'_iv\ar[rr]_-{s'_js'_i\eta'_jv}\ar[d]_-{s'_j\xi_i}
\ar@/^2pc/[rrrr]^-{\lambda'_{i,j}v}
\ar@{}[rrd]|-{\eqref{eq:xi_ij}}&&
s'_js'_is'_jv\ar[r]_-{\eta'_is'_js'_is'_jv}\ar[d]^-{s'_j \zeta_{ij}}&
(s'_is'_j)^2v\ar[r]_-{\mu'_{ij}v}\ar[d]^-{s'_is'_j \zeta_{ij}}&
s'_is'_jv\ar[dddd]_-{ \zeta_{ij}}\ar @/^1pc/[rd]^-{s'_i\xi_j}\\
s'_jvs_i\ar[r]^-{s'_jvs_i\eta_j}\ar[d]_-{\xi_js_i}&
s'_jvs_is_j\ar[r]^-{s'_jv\overline\lambda_{ij}}&
s'_jvs_is_j\ar[r]^-{\eta'_is'_jvs_is_j}\ar[d]^-{\xi_js_is_j}
\ar@{}[rdd]|-{\eqref{eq:xi_ij}}&
s'_is'_jvs_is_j\ar[dd]^-{ \zeta_{ij}s_is_j}&&
s'_ivs_j\ar[dd]^-{\xi_is_j}\\
vs_js_i\ar[r]^-{vs_js_i\eta_j}\ar@{=}[dd]&
vs_js_is_j\ar[r]^-{vs_j\overline \lambda_{ij}}\ar[d]^-{v\eta_is_js_is_j}&
vs_js_is_j\ar[d]^-{v\eta_is_js_is_j}\\
&v(s_is_j)^2\ar[r]^-{vs_is_j \overline \lambda_{ij}}\ar[d]^-{v\mu_{ij}}&
v(s_is_j)^2\ar[r]^-{v\overline \lambda_{ij}s_is_j}&
v(s_is_j)^2\ar[rd]^-{v\mu_{ij}}&&
vs_is_j\ar@/^1pc/[ld]^-{v\overline \lambda_{ij}}\\
vs_js_i\ar[r]_-{v\lambda_{i,j}}&
vs_js_i\ar@{=}[rrr]&&&
vs_is_j}
$$
That is, $(v,\xi_i,\xi_j)$ is a 1-cell in $\Wdl(\ocK)$. Thus the collection
$\{\xi_i:s'_iv\to vs_i\}_{0\leq i \leq n}$ is a 1-cell in
$\Wdl^{(n)}(\ocK)$. The same reasoning as in \eqref{eq:zeta_wreath} shows that
its image under the 2-functor in the claim is the 1-cell $\{\zeta_{\underline 
p}:s_{\underline p}v\to vs_{\underline p}\}_{{\underline p}\in \{0,1\}^{\times(n+1)}}$ 
in $\Mnd(\ocK)^{\two^{n+1}}$ that we started with. Fullness on the 2-cells is
evident.    
\end{proof}

\begin{aussage}{\bf If idempotent 2-cells split.} \label{as:id_split}
Let $\cK$ be a 2-category in which idempotent 2-cells split;
i.e. biequivalent to $\ocK$, and consider the 2-functor
\begin{equation}\label{eq:ktobar}
\xymatrix{
\Wdl^{(n)}(\cK)\ar[r]^-\simeq&
\Wdl^{(n)}(\ocK)\ar[rr]^-{\mathrm{Paragraph~}{\scriptsize\ref{as:embed}}}&&
\Mnd(\ocK)^{\two^{n+1}}\ .}
\end{equation}
By Paragraph \ref{as:embed}, it takes a 1-cell $\{\xi_i:s'_iv\to vs_i\}_{0\leq
i\leq n}$ to an $n+2$-cube in $\Mnd(\ocK)$, with faces of the form in the
first diagram in  
\begin{equation}\label{eq:xi_square}
\xymatrix{
(A,(s_{\underline p+i},\overline\lambda_{\underline p+i}))
\ar[rr]^-{((v,v),\xi_{\underline p+i})}
\ar[d]_{((A,A),\varphi_{\underline p}^i)}&&
(A',(s'_{\underline p+i},\overline{\lambda'}_{\underline p+i})) 
\ar[d]^{((A',A'),{\varphi'}_{\underline p}^i)}&&
(A,z_{\underline p+i})\ar[r]^-{(v,\zeta_{\underline p+i})}
\ar[d]_-{(A,\psi_{\underline p}^i)}&
(A',z'_{\underline p+i})\ar[d]^-{(A',{\psi'_{\underline p}}^i)}\\
(A,(s_{\underline p},\overline\lambda_{\underline p}))
\ar[rr]_-{((v,v),\xi_{\underline p})}&&
(A',(s'_{\underline p},\overline{\lambda'}_{\underline p}))&&
(A,z_{\underline p})\ar[r]_-{(v,\zeta_{\underline p})}&
(A',z'_{\underline p}).}
\end{equation}
Let us choose a biequivalence pseudofunctor $\ocK\to \cK$ adopting the
convention that we split any identity 2-cell trivially; i.e. via identity
2-cells. Then the induced biequivalence $\Mnd(\ocK)\to \Mnd(\cK)$ takes
the first square in \eqref{eq:xi_square} to a {\em commutative} square in
$\Mnd(\cK)$ of the form in the second diagram.
Mapping $\{\xi_i:s'_iv\to vs_i\}_{0\leq i\leq n}$ to the $n+2$-cube in
$\Mnd(\cK)$ formed by these faces; and mapping a 2-cell $\omega$ to the 2-cell
in $\Mnd(\cK)^{\two^{n+1}}$ whose value at each $\underline p\in
\{0,1\}^{\times (n+1)}$ is given by $\omega$; we obtain a fully faithful
2-functor $\Wdl^{(n)}(\cK)\to \Mnd(\cK)^{\two^{n+1}}$. Its composition with
the biequivalence $\Mnd(\cK)^{\two^{n+1}}\stackrel \simeq \to
\Mnd(\ocK)^{\two^{n+1}}$ is 2-naturally isomorphic to \eqref{eq:ktobar}. 
\end{aussage}

\begin{aussage}{\bf If Eilenberg-Moore objects exist.}
Recall (from \cite{street}) that a 2-category $\cK$ is said to admit
Eilenberg-Moore objects provided that the evident inclusion $I:\cK\to
\Mnd(\cK)$ possesses a right 2-adjoint $J$. Whenever $J$ exists, it induces a
fully faithful 2-functor $\Mnd(\cK)\to \cK^\two$ as follows. It takes a monad
$(A,t)$ to the 1-cell part of the counit of the 2-adjunction $I\dashv J$
evaluated at $(A,t)$; i.e. the so-called ``forgetful morphism'' $J(A,t)\to A$
in $\cK$. (The terminology certainly comes from its form in
$\cK=\mathsf{Cat}$.) It takes a 1-cell $(v,\psi)$ to the pair $(v,J(v,\psi))$
and it takes a 2-cell $\omega$ to the pair $(\omega,J\omega)$. 
\end{aussage}

\begin{corollary}
Let $\cK$ be a 2-category in which idempotent 2-cells split and which admits
Eilenberg-Moore objects. Then composing the fully faithful 2-functor
$\Wdl^{(n)}(\cK)\to \Mnd(\cK)^{\two^{n+1}}$ in Paragraph \ref{as:id_split}
with $J^{\two^{n+1}}:\Mnd(\cK)^{\two^{n+1}} \to \cK^{\two^{n+1}}$, we obtain a
fully faithful embedding.  
\end{corollary}

\section{The $n$-ary factorization problem}
\label{sec:factor}

The aim of this section is to find sufficient and necessary conditions on a
demimonad (i.e. a monad in the local idempotent closure of a 2-category) to
be isomorphic to a weak wreath product of $n$ demimonads. 
Some facts about the $n=2$ case are recalled in Paragraph \ref{as:bin_wfac}.

In the next theorem we shall use the notation introduced after Paragraph
\ref{as:K^2n}. 

\begin{theorem}
For any demimonad $(A,s)$ in an arbitrary 2-category $\cK$, the following
assertions are equivalent.  
\begin{itemize}
\item[{(i)}] There is an object $\{\lambda_{i,j}:s_js_i \to s_i s_j\}_{1\leq
  i<j\leq n}$ of $\Wdl^{(n-1)}(\ocK)$ such that the corresponding $n$-ary weak
  wreath product (i.e. its image under the 2-functor in Theorem
\ref{thm:wreath_assoc}) is isomorphic to $(A,s)$.
\item[{(ii)}] There is an $n$ dimensional cube whose 2-faces are commutative
  squares of monad morphisms in $\ocK$ of the form $(A,\varphi^i_{\underline p}):
(A,s_{{\underline p}+i}) \to (A,s_{\underline p})$, such that the
  following hold. For ${\underline p}<{\underline q}\in \{0,1\}^{\times n}$,
  denote by ${\varphi_{\underline p}}^{\underline q}$ and by
  ${\varphi^{\underline p}}_{\underline q}$ the (unique) morphisms composed
  along any path to ${\underline p}+{\underline q}$ from ${\underline p}$ and
  from ${\underline q}$, respectively. Then
\begin{itemize}
\item[{(a)}] $(A,s_{\underline 0})$ is the trivial monad $(A,A)$ and
  $(A,s_{\underline 1})$ is isomorphic to $(A,s)$.
\item[{(b)}] For all ${\underline p}<{\underline q}\in \{0,1\}^{\times n}$, the
  2-cell
$$
\pi_{{\underline p},{\underline q}}:=\big(\xymatrix{
s_{\underline p}s_{\underline q}
\ar[rr]^-{{\varphi_{\underline p}}^{\underline q} \,\,
         {\varphi^{\underline p}}_{\underline q}}&&
s_{{\underline p}+{\underline q}}s_{{\underline p}+{\underline q}}
\ar[r]^-{\mu_{{\underline p}+{\underline q}}}&
s_{{\underline p}+{\underline q}}}\big)
$$
possesses an $s_{\underline p}$-$s_{\underline q}$ bimodule section 
$\iota_{{\underline p},{\underline q}}$ in $\ocK$.
\item[{(c)}]  For all ${\underline p}<{\underline q}<{\underline r}\in
  \{0,1\}^{\times n}$, the morphisms $\iota$ in part (b) render commutative
  the following diagrams. 
\begin{equation}\label{eq:zeta_YB}
\xymatrix{
s_{{\underline p}+{\underline q}+{\underline r}}
\ar[r]^-{\iota_{{\underline p},{\underline q}+{\underline r}}}
\ar[d]_-{\iota_{{\underline p}+{\underline q},{\underline r}}}&
s_{\underline p}s_{{\underline q}+{\underline r}}
\ar[d]^-{s_{\underline p}\iota_{{\underline q},{\underline r}}}\\
s_{{\underline p}+{\underline q}}s_{\underline r}
\ar[r]_-{\iota_{{\underline p},{\underline q}}s_{\underline r}}&
s_{\underline p}s_{\underline q}s_{\underline r}}
\end{equation}
\begin{equation}\label{eq:prod}
\xymatrix @R=15pt{
s_{{\underline q}+{\underline r}}s_{\underline p}
\ar[r]^-{\iota_{{\underline q},{\underline r}}s_{\underline p}}
\ar[dd]_-{{\varphi^{\underline p}}_{{\underline q}+{\underline r}}\,\,
{\varphi_{\underline p}}^{{\underline q}+{\underline r}}}&
s_{\underline q}s_{\underline r}s_{\underline p}
\ar[d]^-{s_{\underline q} {\varphi^{\underline p}}_{\underline r}\,\,     
                       {\varphi_{\underline p}}^{\underline r}}&
s_{\underline r}s_{{\underline p}+{\underline q}}
\ar[r]^-{s_{\underline r}\iota_{{\underline p},{\underline q}}}
\ar[dd]_-{{\varphi^{{\underline p}+{\underline q}}}_{\underline r}\,\,
          {\varphi_{{\underline p}+{\underline q}}}^{\underline r}}&
s_{\underline r}s_{\underline p}s_{\underline q}
\ar[d]^-{{\varphi^{\underline p}}_{\underline r}\,\,     
         {\varphi_{\underline p}}^{\underline r}s_{\underline q}}\\
&s_{\underline q}s_{{\underline p}+{\underline r}}
              s_{{\underline p}+{\underline r}}
\ar[d]^-{s_{\underline q}\mu_{{\underline p}+{\underline r}}}&&
s_{{\underline p}+{\underline r}}
s_{{\underline p}+{\underline r}}s_{\underline q}
\ar[d]^-{\mu_{{\underline p}+{\underline r}}s_{\underline q}}\\
s_{{\underline p}+{\underline q}+{\underline r}}
s_{{\underline p}+{\underline q}+{\underline r}}
\ar[dd]_-{\mu_{{\underline p}+{\underline q}+{\underline r}}}&
s_{\underline q}s_{{\underline p}+{\underline r}}
\ar[d]^-{s_{\underline q}\iota_{{\underline p},{\underline r}}}&
s_{{\underline p}+{\underline q}+{\underline r}}
s_{{\underline p}+{\underline q}+{\underline r}}
\ar[dd]_-{\mu_{{\underline p}+{\underline q}+{\underline r}}}&
s_{{\underline p}+{\underline r}}s_{\underline q}
\ar[d]^-{\iota_{{\underline p},{\underline r}}s_{\underline q}}\\
&s_{\underline q}s_{\underline p}s_{\underline r}
\ar[d]^-{{\varphi^{\underline p}}_{\underline q}\,\,
         {\varphi_{\underline p}}^{\underline q}s_{\underline r}}&&
s_{\underline p}s_{\underline r}s_{\underline q}
\ar[d]^-{s_{\underline p}{\varphi^{\underline q}}_{\underline r}\,\,
                      {\varphi_{\underline q}}^{\underline r}}\\
s_{{\underline p}+{\underline q}+{\underline r}}
\ar[dd]_-{\iota_{{\underline p},{\underline q}+{\underline r}}}&
s_{{\underline p}+{\underline q}}
s_{{\underline p}+{\underline q}}s_{\underline r}
\ar[d]^-{\mu_{{\underline p}+{\underline q}}s_{\underline r}}&
s_{{\underline p}+{\underline q}+{\underline r}}
\ar[dd]_-{\iota_{{\underline p}+{\underline q},{\underline r}}}&
s_{\underline p}s_{{\underline q}+{\underline r}}s_{{\underline q}+{\underline r}}
\ar[d]^-{s_{\underline p}\mu_{{\underline q}+{\underline r}}}\\
&s_{{\underline p}+{\underline q}}s_{\underline r}
\ar[d]^-{\iota_{{\underline p},{\underline q}}s_{\underline r}}&&
s_{\underline p}s_{{\underline q}+{\underline r}}
\ar[d]^-{s_{\underline p}\iota_{{\underline q},{\underline r}}}\\
s_{\underline p}s_{{\underline q}+{\underline r}}
\ar[r]_-{s_{\underline p}\iota_{{\underline q},{\underline r}}}&
s_{\underline p}s_{\underline q}s_{\underline r}&
s_{{\underline p}+{\underline q}}s_{\underline r}
\ar[r]_-{\iota_{{\underline p},{\underline q}}s_{\underline r}}&
s_{\underline p}s_{\underline q}s_{\underline r}}
\end{equation}
\end{itemize}
\end{itemize}
\end{theorem}

\begin{proof}
\underline{(i)$\Rightarrow$(ii).}
The 2-functor in Paragraph \ref{as:embed} takes $\{\lambda_{i,j}:s_js_i \to
s_i s_j\}_{1\leq i<j\leq n}$  to a commutative $n$-cube in $\Mnd(\ocK)$ with
edges $(A,\varphi_{\underline p}^i):(A,s_{{\underline p}+i})\to
(A,s_{\underline p})$ of the form in Lemma \ref{lem:A}. In this cube
$(A,s_{\underline 0})$ is the trivial monad $(A,A)$ and $(A,s_{\underline 1})$
is the $n$-ary weak wreath product which is isomorphic to $(A,s)$ by
assumption. Thus property (a) holds. By construction of the 2-functor in
Paragraph \ref{as:embed}, $(A,s_{{\underline p}+{\underline q}})$ is the weak
wreath product of $(A,s_{\underline p})$ and $(A,s_{\underline q})$, for all
${\underline p}< {\underline q}\in \{0,1\}^{\times n}$. Hence the 2-cell
$\pi_{{\underline p}, {\underline q}}$ in part (b) possesses a bilinear
section $\iota_{{\underline p}, {\underline q}}$ by Paragraph
\ref{as:bin_wfac}. It remains to show that the diagrams in part (c) commute.  

The monic 2-cell
$\iota_{{\underline p}, {\underline q}}$
is given by 
$\overline \lambda_{{\underline p} +{\underline q}}:
(s_{\underline p} s_{\underline q},
\overline \lambda_{{\underline p}+{\underline q}})
\to (s_{\underline p} s_{\underline q},
\overline s_{\underline p} \overline s_{\underline q})$; 
$\iota_{{\underline p}, {\underline q}+{\underline r}}$
is equal to
$\overline \lambda_{{\underline p}+{\underline q}+{\underline r}}:
(s_{\underline p} s_{\underline q} s_{\underline r},
\overline \lambda_{{\underline p}+{\underline q}+{\underline r}})
\to (s_{\underline p} s_{\underline q} s_{\underline r}, 
\overline s_{\underline p}\overline \lambda_{{\underline q}+
{\underline r}} )$ and so on. Hence \eqref{eq:zeta_YB} takes the form
$$
\xymatrix @C=60pt {
(s_{\underline p} s_{\underline q} s_{\underline r},
\overline \lambda_{{\underline p}+{\underline q}+{\underline r}})
\ar[r]^-{\overline \lambda_{{\underline p}+{\underline q}+{\underline r}}}
\ar[d]_-{\overline \lambda_{{\underline p}+{\underline q}+{\underline r}}} &
(s_{\underline p} s_{\underline q} s_{\underline r},
\overline s_{\underline p}\overline \lambda_{{\underline q}+{\underline r}})
\ar[d]^-{s_{\underline p}\overline \lambda_{{\underline q}+{\underline r}}}\\
(s_{\underline p} s_{\underline q} s_{\underline r},
\overline \lambda_{{\underline p}+{\underline q}}s_{\underline r})
\ar[r]_-{\overline \lambda_{{\underline p}+{\underline q}}
\overline s_{\underline r}}&
(s_{\underline p} s_{\underline q} s_{\underline r},
s_{\underline p} s_{\underline q} s_{\underline r})}
$$
which commutes by \eqref{eq:rarrownorm} and \eqref{eq:larrownorm}.
In the vertical paths of the diagrams in \eqref{eq:prod}, note the occurrence
of the weak distributive laws $\lambda_{{\underline p},{\underline q}}$,
$\lambda_{{\underline p},{\underline q}+{\underline r}}$, etc. Thus the first
diagram in \eqref{eq:prod} takes the form 
$$
\xymatrix @C=40pt {
(s_{\underline q} s_{\underline r}s_{\underline p},
\overline \lambda_{{\underline q}+{\underline r}}\overline s_{\underline p})
\ar[d]_-{\overline \lambda_{{\underline q}+{\underline r}}s_{\underline p}}
\ar[rr]^-{\lambda_{{\underline p},{\underline q}+{\underline r}}}&&
(s_{\underline p} s_{\underline q} s_{\underline r},
\overline s_{\underline p}\overline \lambda_{{\underline q}+{\underline r}})
\ar[d]^-{s_{\underline p}\overline \lambda_{{\underline q}+{\underline r}}}\\
(s_{\underline q} s_{\underline r}s_{\underline p},
\overline s_{\underline q}\overline s_{\underline r}\overline s_{\underline p})
\ar[r]_-{s_{\underline q}\lambda_{{\underline p},{\underline r}}}&
(s_{\underline q} s_{\underline p}s_{\underline r},
\overline s_{\underline q}\overline s_{\underline p}\overline s_{\underline r})
\ar[r]_-{\lambda_{{\underline p},{\underline q}}s_{\underline r}}&
(s_{\underline p} s_{\underline q} s_{\underline r},
\overline s_{\underline p}\overline s_{\underline q}\overline s_{\underline r})}
$$
which is evidently commutative in view of \eqref{eq:Ck0}.
Commutativity of the second diagram in \eqref{eq:prod}
follows symmetrically.  

\underline{(ii)$\Rightarrow$(i).}
By assumption, for all $1\leq i<j\leq n$, the 2-cell
$$
\pi_{i,j}:=\big(\xymatrix{
s_is_j\ar[rr]^-{{\varphi_i}^j\,\,{\varphi^i}_j}&&
s_{ij}s_{ij}\ar[r]^-{\mu_{ij}}&
s_{ij}}\big)
$$
possesses an $s_i$-$s_j$ bimodule section $\iota_{i,j}$. Hence by Paragraph
\ref{as:bin_wfac}, there is a weak distributive law $\lambda_{i,j}:=
\iota_{i,j}.\mu_{ij}.{\varphi^i}_j\,\,{\varphi_i}^j:s_js_i\to s_i s_j$ in
$\ocK$ such that the corresponding weak wreath product is isomorphic to
$s_{ij}$. Let us prove that the collection $\{\lambda_{i,j}:s_js_i \to s_i
s_j\}_{1\leq i<j\leq n}$ obeys the Yang-Baxter conditions, for all $1\leq i<j<
k \leq n$. This follows by commutativity of the following diagram. 
$$
\scalebox{.95}{
\xymatrix @R=15pt {
s_ks_js_i \ar[rr]^-{s_k {\varphi^i}_j\,\,{\varphi_i}^j}
\ar[dd]_-{{\varphi^j}_k\,\,{\varphi_j}^ks_i}&&
s_ks_{ij}s_{ij}\ar[r]^-{s_k\mu_{ij}}
\ar[dd]^-{{\varphi^{ij}}_k\,\,{\varphi_{ij}}^k\,\,{\varphi_{ij}}^k}&
s_ks_{ij}\ar[rr]^-{s_k\iota_{i,j}}
\ar[dd]^-{{\varphi^{ij}}_k\,\,{\varphi_{ij}}^k}
\ar@{}[rrdddd]|-{\eqref{eq:prod}}&&
s_ks_is_j 
\ar[d]^-{{\varphi^i}_k\,\,{\varphi_i}^k s_j}\\
&&&&&
s_{ik}s_{ik}s_j\ar[d]^-{\mu_{ik}s_j}\\
s_{jk}s_{jk}s_i
\ar[rr]^-{{\varphi^i}_{jk}\,\,{\varphi^i}_{jk}\,\,{\varphi_i}^{jk}}
\ar[dd]_-{\mu_{jk}s_i}&&
s_{ijk}s_{ijk}s_{ijk}\ar[r]^-{s_{ijk}\mu_{ijk}}
\ar[dd]^-{\mu_{ijk}s_{ijk}}&
s_{ijk}s_{ijk}\ar[dd]^-{\mu_{ijk}}&&
s_{ik}s_j\ar[d]^-{\iota_{i,k}s_j}\\
&&&&&
s_is_ks_j\ar[d]^-{s_i {\varphi^j}_k\,\,{\varphi_j}^k}\\
s_{jk}s_i \ar[rr]^-{{\varphi^i}_{jk}\,\,{\varphi_i}^{jk}}
\ar[d]_-{\iota_{j,k}s_i}
\ar@{}[rrrdd]|-{\eqref{eq:prod}}&&
s_{ijk}s_{ijk}\ar[r]^-{\mu_{ijk}}&
s_{ijk}\ar[r]^-{\iota_{ij,k}}\ar[d]_-{\iota_{i,jk}}
\ar@{}[rrdd]|-{\eqref{eq:zeta_YB}}&
s_{ij}s_k\ar[rdd]^(.3){\iota_{i,j}s_k}&
s_i s_{jk} s_{jk} \ar[d]^-{s_i \mu_{jk}}\\
s_js_ks_i \ar[d]_-{s_j {\varphi^i}_k\,\,{\varphi_i}^k}&&&
s_i s_{jk}\ar[rrd]_(.3){s_i \iota_{j,k}}&&
s_i s_{jk}\ar[d]^-{s_i \iota_{j,k}}\\
s_j s_{ik}s_{ik}\ar[r]_-{s_j \mu_{ik}}&
s_j s_{ik}\ar[r]_-{s_j \iota_{i,k}}&
s_js_is_k \ar[r]_-{{\varphi^i}_j\,\,{\varphi_i}^js_k}&
s_{ij}s_{ij}s_k\ar[r]_-{\mu_{ij}s_k}&
s_{ij}s_k\ar[r]_-{\iota_{i,j}s_k}&
s_is_js_k}}
$$
Thus $\{\lambda_{i,j}:s_js_i \to s_i s_j\}_{1\leq i<j\leq
n}$ is an object of $\Wdl^{(n-1)}(\ocK)$. It remains to show that the
corresponding $n$-ary weak wreath product is isomorphic to $(A,s_{\underline
1})\cong (A,s)$.

As observed above, $(A,s_{12})$ is isomorphic to the weak wreath product of
$(A,s_1)$ and $(A,s_2)$ with respect to the weak distributive law
$\lambda_{1,2}:= \iota_{1,2}.\mu_{12}.{\varphi^1}_2\,\,{\varphi_1}^2:s_2s_1\to
s_1 s_2$. Similarly, $(A,s_{123})$ is isomorphic to the weak wreath product of
$(A,s_{12})$ and $(A,s_3)$ with respect to $\lambda_{12,3}:=
\iota_{12,3}.\mu_{123}.{\varphi^{12}}_3\,\,{\varphi_{12}}^3:s_3s_{12}\to
s_{12} s_3$. Moreover, precomposing both paths around the second
diagram in \eqref{eq:prod} (for $\underline p=1$, $\underline q=2$ and
$\underline r=3$) by $s_3\pi_{1,2}=s_3\mu_{12}.s_3{\varphi_1}^2{\varphi^1}_2$, 
we obtain that the weak distributive law $\lambda_{12,3}:(s_3s_{12},\overline
s_3\overline s_{12})\to (s_{12}s_3,\overline s_{12}\overline s_3)$ differs by
the isomorphisms  
$\xymatrix{
\iota_{1,2}:(s_{12},\overline s_{12})\ar@<2pt>[r]&
(s_1s_2,\overline \lambda_{12}):\pi_{1,2}\ar@<2pt>[l]}$ 
from the image
$s_1\lambda_{2,3}.\lambda_{1,3}s_2.s_3\overline \lambda_{12}:$ 
$(s_3s_1s_2, s_3\overline \lambda_{12})\to (s_1s_2s_3,\overline
\lambda_{12}s_3)$  
of $\{\lambda_{i,j}:s_js_i\to s_is_j\}_{1\leq i<j\leq 3}$ under
the 2-functor $C_1:\Wdl^{(2)}(\ocK)\to \Wdl^{(1)}(\ocK)$ in Theorem
\ref{thm:n_w-wreath}.
Hence $(A,s_{123})$ is isomorphic to the ternary weak wreath product of
$\{(A,s_i)\}_{1\leq i \leq 3}$. Iterating this reasoning we conclude that
$(A,s_{\underline 1})$ is isomorphic to the $n$-ary weak wreath product of
$\{(A,s_i)\}_{1\leq i \leq n}$.  
\end{proof}


\begin{thebibliography}{Bibliography}{}

\bibitem{AAatAl}
J.N. Alonso \'Alvarez, J.M. Fern\'andez Vilaboa and R. Gonz\'alez
Rodr\'{\i}guez. 
{\em Weak Hopf algebras and weak Yang-Baxter operators}. 
J. Algebra 320 (2008), no. 5, 2101-2143.

\bibitem{Beck} 
J. Beck, 
{\em Distributive laws}.
in: {\em Seminar on triples  and categorical homology theory}, 
B. Eckmann (ed.), Springer LNM 80, 119-140 (1969).    

\bibitem{Bohm:PhD}
G. B\"ohm. 
{\em Weak Hopf algebras and their application to spin models}.
PhD thesis, Budapest, 1997.

\bibitem{Boh_wtm} 
G. B\"ohm, {\em The weak theory of monads}.
Adv. in Math. 225 (2010), 1-32. 

\bibitem{BoCaJa:wba}
G. B\"ohm, S. Caenepeel and K. Janssen.
{\em Weak bialgebras and monoidal categories}. 
Comm. Algebra 39 (2011), no. 12 (special volume dedicated to Mia Cohen),
4584-4607.  

\bibitem{BohmGomTor}
G. B\"ohm and J. G\'omez-Torrecillas,
{\em Weak factorization of algebras}.
Preprint available at \href{http://arxiv.org/abs/1108.5957}{arXiv:1108.5957v1}.

\bibitem{BoLaSt_wdlaw}
G. B\"ohm, S. Lack and R. Street.
{\em On the 2-categories of weak distributive laws}.
Comm. Algebra 39 (2011), no. 12 (special volume dedicated to Mia Cohen),
4567-4583. 

\bibitem{BoLaSt_weak}
G. B\"ohm, S. Lack and R. Street.
{\em Idempotent splittings, colimit completion, and weak aspects of the theory
  of monads}. 
J. Pure Appl. Algebra, 216 (2012) no. 2, 385-403. 

\bibitem{BNSz}
G. B\"ohm, F. Nill, and K. Szlach\'anyi.
{\em Weak Hopf algebras. I. Integral theory and $C^*$-structure}.
J. Algebra 221 (1999), no. 2, 385-438.

\bibitem{Stef&Erwin} 
S. Caenepeel and E. De Groot, 
{\em Modules over weak entwining structures}.
[in:] {\em New Trends in Hopf Algebra Theory}, 
N. Andruskiewitsch, W.R. Ferrer Santos and H-J. Schneider (eds.),
Contemp. Math. 267 (2000), 4701-4735.

\bibitem{Cheng:dlaws}
E. Cheng, 
{\em Iterated distributive laws}.
Mathematical Proceedings of the Cambridge Philosophical Society 150 (2011),
459-487.  

\bibitem{Nill}
F. Nill.
{\em Weak bialgebras}. 
Preprint available at 
\href{http://arxiv.org/abs/math/9805104}{arXiv:math/9805104v1}.

\bibitem{NilSzl:Hopf_spin}
F. Nill and K. Szlach\'anyi.
{\em Quantum chains of Hopf algebras with quantum double cosymmetry}.
Commun. Math. Phys. 187 (1997), 159-200.

\bibitem{NilSzlWie:Jones}
F. Nill,  K. Szlach\'anyi and H.-W. Wiesbrock.
{\em  Weak Hopf algebras and reducible Jones inclusions of depth 2. I: From
  crossed products to Jones towers}. 
Preprint available at
\href{http://arxiv.org/abs/math/9806130}{arXiv:math/9806130v1}. 

\bibitem{PaSt}
C. Pastro and R. Street. 
{\em Weak Hopf monoids in braided monoidal categories}. 
Algebra Number Theory 3 (2009), no. 2, 149-207.

\bibitem{street}
R.H. Street,
{\em The formal theory of monads}.
J. Pure Appl. Algebra 2 (1972), 149-168.

\bibitem{Str:weak_dl} 
R.H. Street, 
{\em Weak distributive laws}.
Theory and Applications of Categories 22 (2009), 313-320.

\bibitem{SzlVec:G-spin}
K. Szlach\'anyi and P. Vecserny\'es, 
{\em Quantum symmetry and braid group statistics in G-spin models}.
Commun. Math. Phys. 156 (1993), 127-168.

\end{thebibliography}
\end{document}